\documentclass[reqno,10pt]{amsart}

\usepackage[utf8]{inputenc}
\usepackage[T1]{fontenc}

\usepackage{amsmath}
\usepackage{amssymb}
\usepackage{amsthm}
\usepackage[english]{babel}
\usepackage{bbm}
\usepackage{booktabs} 
\usepackage{cancel}
\usepackage{enumitem}
\usepackage{flafter} 
\usepackage{graphicx}
\usepackage[margin=2.54cm]{geometry}
\usepackage{mathtools}
\usepackage{microtype} 
\usepackage{multirow} 
\usepackage[square, numbers]{natbib}
\usepackage{pgfplots}
\usepackage{placeins} 
\usepackage{setspace}
\usepackage{siunitx} 
\usepackage{subcaption} 
\usepackage{tikz}
\usepackage{xspace} 

\usepackage[hidelinks]{hyperref}

\let\citet\citep

\usepgfplotslibrary{groupplots} 
\usepgfplotslibrary{statistics} 

\pgfplotsset{
  every axis/.append style={
    label style={font=\small},
    ticklabel style={font=\footnotesize},
    legend style={font=\small},
    title style={font=\small}
  }
}

\theoremstyle{plain}
\usepackage{amsthm}
\newtheorem{theorem}{Theorem}
\newtheorem{proposition}[theorem]{Proposition}

\newtheorem{lemma}[theorem]{Lemma}

\theoremstyle{definition}

\theoremstyle{remark}

\numberwithin{equation}{section}
\numberwithin{theorem}{section}

\theoremstyle{plain}

\newtheorem*{assumptionS*}{Assumption S}

\newtheorem*{assumptionT*}{Assumption T}

\newtheorem*{assumptionU*}{Assumption U}

\newtheorem*{assumptionR*}{Assumption R}

\newtheorem*{assumptionP*}{Assumption P}

\newtheorem*{assumptionC1*}{Assumption C1}

\newtheorem*{assumptionC2*}{Assumption C2}

\newtheorem*{modelA*}{Model A}

\newtheorem*{modelB*}{Model B}

\def\bs{\boldsymbol}
\def\mrm{\mathrm}
\def\msf{\mathsf}
\def\mcl{\mathcal}
\def\mbb{\mathbb}
\def\mbf{\mathbf}

\let\P\undefined
\let\exp\undefined
\let\log\undefined
\def\d{\mathop{}\!\mrm d}

\def\deg{\mathop{}\!\msf{deg}}
\DeclarePairedDelimiterXPP\P[1]{\mathop{}\!\mbb{P}}{(}{)}{}{#1}
\DeclarePairedDelimiterXPP\E[1]{\mathop{}\!\mbb{E}}{[}{]}{}{#1}
\DeclarePairedDelimiterXPP\1[1]{\mathop{}\!\mathbbm{1}}{\{}{\}}{}{#1}
\DeclarePairedDelimiterXPP\Poi[1]{\mathop{}\!\msf{Poi}}{(}{)}{}{#1}
\DeclarePairedDelimiterXPP\Var[1]{\mathop{}\!\msf{Var}}{(}{)}{}{#1}
\DeclarePairedDelimiterXPP\Cov[1]{\mathop{}\!\msf{Cov}}{(}{)}{}{#1}
\DeclarePairedDelimiterXPP\exp[1]{\mathop{}\!\msf{exp}}{(}{)}{}{#1}
\DeclarePairedDelimiterXPP\log[1]{\mathop{}\!\msf{log}}{(}{)}{}{#1}

\def\f{\frac}
\def\e{\mrm{e}}
\def\w{\wedge}
\def\any{\;\cdot\;}
\def\es{\varnothing}
\def\su{\subseteq}
\def\co{\colon}
\def\coeq{\coloneq}
\def\sm{\setminus}
\def\ff{\infty}
\def\tff{\uparrow\ff}
\def\given{\nonscript\:\vert\nonscript\:\mathopen{}}
\def\biggiven{\nonscript\:\big\vert\nonscript\:\mathopen{}}

\def\ti{\times}

\def\NN{\mcl N}
\def\PP{\mcl P}
\def\SS{\mcl S}

\def\le{\leqslant}
\def\ge{\geqslant}

\DeclarePairedDelimiter{\abs}{\lvert}{\rvert}

\DeclarePairedDelimiter{\ceil}{\lceil}{\rceil}
\DeclarePairedDelimiter{\set}{\{}{\}}
\def\tf{\tfrac}

\def\R{\mbb R}

\def\S{\mbb S}

\def\pp{\bs{p}}

\def\a{\alpha}
\def\b{\beta}
\def\g{\gamma}
\def\Ga{\Gamma}
\def\de{\delta}
\def\De{\Delta}

\def\la{\lambda}

\def\r{\rho}
\def\t{\tau}

\def\s{\sigma}

\def\Nl{\mbf{N}_{\msf{loc}}}

\def\Gh{G^{\msf{Dow}}}
\def\Gb{G^{\msf{bip}}}
\def\sna{S_n^\ge}
\def\snb{S_n^\le}
\def\snc{S_n^{(1)}}
\def\snd{S_n^{(2)}}
\def\vpp{\overrightarrow{\bs{p}}}

\def\vPp{\overrightarrow{\bs{P}}}
\def\vPs{\overrightarrow{\bs{P}'}}

\makeatletter

\addto\extrasenglish{}
\addto\extrasenglish{}
\makeatother

\ExplSyntaxOn

\NewDocumentCommand\MakeHatBraceCommand{mm}{
  \NewDocumentCommand#1{O{#2}m}{
    \bool_lazy_and:nnTF
      { \tl_if_single_token_p:n {##2} } { \token_if_macro_p:N ##2 }
      { \__daleif_bh_gen:nno ##1 ##2 }
      { \__daleif_bh_gen:nnn ##1 {##2} }
  }
}

\cs_new_protected:Npn \__daleif_bh_gen:nnn #1#2#3
  {
    \bool_lazy_and:nnTF
      { \tl_if_single_p:n {#3} } { \tl_if_head_is_group_p:n {#3} }
      { \__daleif_bh_gen_process:nno #1 #2 { \use:n #3 } }
      { \__daleif_bh_gen_process:nnn #1 #2 {#3} }
  }
\cs_generate_variant:Nn \__daleif_bh_gen:nnn { nno }
\cs_new_protected:Npn \__daleif_bh_gen_process:nnn #1 #2 #3
  {
    \tl_if_head_eq_meaning:nNTF {#3} #1
      { \tl_set:Nx \l_tmpa_tl { \tl_tail:n {#3} } }
      { \tl_set:Nn \l_tmpa_tl {#3} }
    \tl_if_in:onT { \l_tmpa_tl \q_stop } { #2\q_stop }
      {
        \tl_set:No \l_tmpa_tl { \l_tmpa_tl \q_stop }
        \tl_replace_once:Nnn \l_tmpa_tl { #2\q_stop } {}
      }
    \sp { \smash {#1} \smash[b]{ \l_tmpa_tl } \smash {#2} }
  }
\cs_generate_variant:Nn \__daleif_bh__gen_process:nnn { nno }

\MakeHatBraceCommand\ph{()}

\ExplSyntaxOff

\renewcommand{\binom}[2]{\genfrac{}{}{0pt}{}{#1}{#2}}
\def\Cpp{C\texttt{++}\xspace}
\def\xmin{x_\msf{\scriptstyle min}}

\pgfmathdeclarefunction{normalpdf}{3}{\pgfmathparse{1/(#3*sqrt(2*pi))*exp(-((#1-#2)^2)/(2*#3^2))}}

\author[M. Brun]{Morten Brun$^1$}
\email{morten.brun@uib.no}

\author[C. Hirsch]{Christian Hirsch$^2$}
\email{hirsch@math.au.dk}

\author[P. Juhasz]{Peter Juhasz$^{2\ast}$}
\email{peter.juhasz@math.au.dk}

\author[M. Otto]{Moritz Otto$^3$}
\email{m.f.p.otto@math.leidenuniv.nl}

\address{$^1$Department of Mathematics, University of Bergen, Bergen, Norway}
\address{$^2$Department of Mathematics, Aarhus University, Aarhus, Denmark}
\address{$^3$Mathematical Institute, Leiden University, Leiden, Netherlands}
\address{$^\ast$Corresponding author}

\begin{document}

\title{Random Connection Hypergraphs}


\begin{abstract}
    In this paper, we introduce a novel model for random hypergraphs based on weighted random connection models.
    In accordance with the standard theory for hypergraphs, this model is constructed from a bipartite graph.
    In our stochastic model, both vertex sets of this bipartite graph form marked Poisson point processes, and the connection radius is inversely proportional to a product of suitable powers of the marks.
    Hence, our model is a common generalization of weighted random connection models and AB random geometric graphs.
    For this hypergraph model, we investigate the limit theory of various graph-theoretic and topological characteristics, including higher-order degree distributions, Betti numbers of the associated Dowker complex, and simplex counts.
    In particular, for the latter quantity, we identify regimes of convergence to normal and to stable distribution depending on the heavy-tailedness of the weight distribution.
    We conclude our investigation with a simulation study and an application to the collaboration network extracted from the arXiv dataset.
\end{abstract}

\maketitle
\subsection*{Keywords} Poisson point process; random graph; hypergraph; simplicial complex; central limit theorem

\section{Introduction}\label{sec:int}

%
%
Since the seminal work of Barab\'asi and Albert~\citep{barabasi}, the study of complex networks has become a very active field of research.
The primary reason for this is that many real-world systems can be represented as networks, where nodes represent the elements of the system and links represent the interactions between them.
Examples of such systems include the Internet, social networks, and biological networks, among others.
The study of complex networks has led to the development of new tools and methods for analyzing and modeling the structure and dynamics of these systems.

%
%
The popularity of complex network models can be attributed to their ability to capture the main features of real-world networks such as the presence of a scale-free degree distribution and the small-world property, which means that the average distance between any two nodes in the network is relatively short.
Despite these promising results, the most basic complex network models struggle to capture clustering structures that are essential in various application domains.
An elegant approach to implementing these clustering effects involves embedding the network nodes in a suitable ambient space, where nearby nodes have a greater tendency to be connected.
This geometric embedding, therefore, encourages the formation of clusters of nodes in a spatial vicinity.
While there are various ways to endow complex networks with a spatial structure, one of the most elegant and mathematically tractable approaches is the age-dependent random connection model (ADRCM) in \citet{glm2, glm}.
These models enable the combination of a geometric embedding with possibly heavy-tailed vertex weights, which can create a hierarchy of hubs often observed in real data.
This idea also lies at the core of kernel-based networks as considered in \citet{gracar2022recurrence, komjathy}.

%
%
Besides the spatial structure, another important advantage of the ADRCM is that it can be used to model higher-order interactions in complex networks.
For instance, in the context of collaboration networks, scientific papers often involve more than one or two authors, thereby illustrating the need to go beyond binary networks.
One of the challenges in investigating higher-order networks is that their topology can become highly complex, which is why tools from the field of algebraic topology are now being applied in the context of network analysis.
For instance, \citet{siu2023many} considered a higher-order model for preferential attachment and computed the asymptotic growth rate of the expected Betti numbers.
Moreover, in an earlier work, we studied the potential of the ADRCM as a model for arXiv collaboration data from different fields~\citep{hj23}.

%
%
While the ADRCM is an exciting model for scale-free networks with a spatial structure, it is too coarse to capture some key features in collaboration networks.
Indeed, from the ADRCM, we can tell whether two scientists have collaborated, but there is no way to determine how many papers they have written together.
To address this shortcoming, we propose a spatial hypergraph structure, which we call the \emph{random connection hypergraph model (RCHM)}.
This hypergraph, which we think of as a collection of subsets of vertices, will be defined through its incidence graph, i.e., the bipartite graph whose two vertex sets are given by the set of hypergraph vertices and relations, respectively~\citep{berge}.
In our setting, this bipartite graph has vertex sets corresponding to the authors and documents, respectively, with an edge indicating that an author has collaborated on a document.

On the model side, the RCHM removes one of the key drawbacks of the ADRCM, namely that the latter does not allow for capturing the number of documents that a collection of authors has collaborated on.
Despite this improvement, we show that the RCHM remains mathematically tractable.
Indeed, we can recover many of the central higher-order characteristics considered in \citet{hj23} such as higher-order degree distributions, simplex counts, and Betti numbers.
However, the bipartite structure often induces additional complexities that require substantially more involved mathematical machinery compared to the ADRCM.
Surprisingly, in the case of simplex counts, the bipartite structure allows us to go further than the ADRCM.

Hence, the main contributions of this manuscript are the following:
\begin{itemize}
    \item
        We introduce the RCHM, a higher-order model for bipartite networks.
        Note that in the setting of traditional spatial networks, the \emph{AB random geometric graphs} are an intensively studied model~\citep{ab3, ab4, ab1, ab2}.
        However, while bipartite AB random geometric bipartite graphs have previously received substantial attention in telecommunication networks, they are not suited for collaboration networks, as they are not scale-free.
        To the best of our knowledge, the RCHM is the first spatial model applicable to scale-free bipartite networks.
    \item
        While the bipartite structure of the RCHM induces additional complexities, we show that the RCHM is mathematically tractable.
        In fact, surprisingly, in the case of simplex counts, the bipartite structure allows us to go further than the ADRCM.
        We note that the triangle count is of high interest in the context of complex networks since it is intimately related to the clustering coefficients.
        For both geometric and combinatorial network models, the clustering coefficient has recently been the subject of intense research~\citep{litvak, pvh}.
    \item
        We illustrate the usefulness of our asymptotic results in finite sample sizes through an extensive simulation study.
        Finally, we also compare the model to real-world data from the arXiv network.
\end{itemize}

The rest of the manuscript is organized as follows.
In \autoref{sec:mod}, we introduce the RCHM and discuss its main properties.
We also present the main results of the manuscript.
Sections~\ref{sec:palm_distribution}--\ref{sec:sim} are devoted to the proofs of the main results.
Next, in \autoref{sec:sim} we illustrate the applicability of our asymptotic results in finite sample sizes through an extensive simulation study.
Finally, in \autoref{sec:dat}, we compare the model to real-world data from the arXiv network.

\section{Model and main results}\label{sec:mod}

We work on the space $\S \coeq \R \ti [0, 1]$ and interpret~$\R$ as the \emph{location} and $[0, 1]$ as the \emph{mark space}.
Given parameters $0 < \g, \g' < 1$, $\b > 0$, we introduce our notion of connectivity of elements in~$\S$.
For $p \coeq (x, u) \in \S$ and $p' \coeq (z, w) \in \S$ we let
\begin{equation} \label{eq:F} F(p, p') = F(p, p'; \g, \g') \coeq \abs{x - z} u^{\g} w^{\g'} \end{equation}
and say that $p \in \S$ and $p' \in \S$ are \emph{connected} if
\begin{equation} F(p, p'; \g, \g') \le \b. \label{eq:wrchm} \end{equation}
We define the \emph{neighborhood} of a point $p \in \S$ by
\begin{equation} B(p, \b) = B_F(p, \b) \coeq \set[\big]{p' \in \S \co F(p, p'; \g, \g') \le \b}. \label{eq:ball} \end{equation}
Loosely speaking, we can think of~$B$ as a non-metric ball centered at a point $p \in \S$ with radius~$\b$.
For $\De \su \S$, we also put $B_F(\De, \b) \coeq \bigcap_{p \in \De} B_F(p, \b)$ as the joint neighborhood of the points in~$\De$.

Next, we define a random bipartite graph.
Let $\PP, \PP'$ be independent Poisson processes on~$\S$ with intensity measures $\la \abs{\any}, \la' \abs{\any}$, respectively, where $\la, \la' > 0$ are fixed parameters and $\abs{\any}$ denotes Lebesgue measure on~$\S$.
Let $\Gb \coeq \Gb(\PP, \PP') \coeq \Gb(\PP, \PP'; \g, \g', \b)$ be the bipartite graph with vertex set $\PP \cup \PP'$, where there is an edge between two points $p \in \PP$ and $p' \in \PP'$ if and only if~$p$ and~$p'$ are connected.
Note that given~$\PP$ and~$\PP'$, the number of edges in~$\Gb$ is increasing in the parameters $\g, \g'$ and~$\b$ and that $\PP' \cap B_F(p, \b)$ is the subset of points in~$\PP'$ connected to~$p$ in~$\Gb$.

In the collaboration network example discussed in \autoref{sec:dat} below, the Poisson processes~$\PP$ and~$\PP'$ can be thought of as sets of authors and documents, respectively.
An edge $(p, p')$ can be interpreted as the relation that author $p \in \PP$ has collaborated on manuscript $p' \in \PP'$.
Considering the product on the right-hand side of~\eqref{eq:F} illustrates that~$\Gb$ is a natural bipartite extension of the age-dependent random connection models from \citet{glm2, glm}.
As will be made precise below, $\g, \g' < 1$ determine the power-law exponent of the degree distribution for vertices in $\PP, \PP'$, respectively.
Once $\g, \g'$ are fixed, the parameter $\b > 0$ can be used to tune the overall number of expected edges in the bipartite graph~$\Gb$.

%
%
Note that it follows from identity~\eqref{eq:F} that the roles of~$\PP$ and~$\PP'$ are symmetric when switching the parameters~$\g$ and~$\g'$.
Hence, when in the following we describe more elaborate graph-theoretic and topological quantities of~$\Gb$ that are defined with reference to~$\PP$, the symmetry between~$\PP$ and~$\PP'$ implies that our results also hold for the corresponding quantities defined in terms of~$\PP'$, when switching the assumptions on~$\g$ and~$\g'$.

%
%
We now introduce the \emph{random connection hypergraph model (RCHM)} $\Gh \coeq \Gh(\PP, \PP')$ through the construction known as \emph{Dowker complex}~\citep{dowker}.
We stress here that~$\Gh$ is a specific type of hypergraph, namely, a simplicial complex.
In the following, we rely on this structure as a simplicial complex to investigate Betti numbers.
More precisely, the set of hyperedges $\Sigma_m$ of~$\Gh$ of cardinality $m + 1$ is given by
\begin{equation} \Sigma_m \coeq \set[\big]{\De_m \su \PP \co \#(\De_m) = m + 1, \PP' \cap B_F(\De_m, \b) \ne \es}, \label{eq:simplices} \end{equation}
where $\#(\any)$ denotes the cardinality of a set.

Following the terminology in topological data analysis, we also refer to the elements of $\Sigma_m$ as \emph{$m$-simplices}.
If there is a unique vertex of an $m$-simplex $\De_m$ with the lowest mark, then we define $c(\De_m) \in \R$ as the location coordinate of this vertex.
Otherwise, let $c(\De_m)$ be the location of the left-most vertex among all vertices with minimal mark.
Although not explicitly stated, the set of $m$-simplices $\Sigma_m$ depends on the point process~$\PP'$.
In our interpretation, $m + 1$ authors form an $m$-simplex if and only if they have coauthored at least one common paper.

%
%
In classical binary networks, the degree distribution of a typical vertex is of central importance.
Due to the translation-invariance of the edge rule~\eqref{eq:wrchm}, such a vertex can be chosen to be of the form $o = (0, U)$ with~$U$ uniform in $[0, 1]$.
For higher-order networks, it is essential to go beyond single vertices and describe properties of typical $m$-simplices.
However, rigorously defining the notion of a typical simplex is a far more involved process.
To address this problem, we will rely on the established theory of Palm calculus as explained in \citet[Chapter~9]{poisBook}.
As a first step, we need to establish that the $m$-simplex intensity is finite.
Given a subset~$A$ of~$\R$, we denote the subset of $m$-simplices centered in~$A$ by
\[ \Sigma^A_m \coeq \set[\big]{\De_m \in \Sigma_m \co c(\De_m) \in A}. \]
\begin{proposition}[Finiteness of the $m$-simplex-intensity]\label{prop:finite_intensity}
    Let~$\PP$ and~$\PP'$ be independent Poisson point processes on~$\S$, let $A \su \R$ be a Borel set with Lebesgue measure one.
    Given $m \ge 0$, $\g < 1$, $\g' < 1/(m + 1)$ and $\b > 0$, the $m$-simplex intensity $\la_m \coeq \E{\# \Sigma^A_m}$ is a non-zero finite value not depending on the choice of~$A$.
\end{proposition}

Now, we define the distribution of the typical $m$-simplex $\De_m\ph{\ast}$.
To do so, we first introduce the following notations.
More precisely,~$f$ denotes an arbitrary nonnegative measurable functional that may depend on a considered $m$-simplex as well as on the point processes~$\PP$ and~$\PP'$.
Then, we define the distribution of the typical $m$-simplex as
\begin{equation} \E[\big]{f(\De_m\ph{\ast}, \PP, \PP')} = \f{1}{\la_m} \E[\Big]{\sum_{\De_m \in \Sigma_m^{A}} f \bigl( \De_m - c(\De_m), \PP - c(\De_m), \PP' - c(\De_m) \bigr)} \label{eq:palm} \end{equation}
where $A \su \R$ is an arbitrary Borel set with $\abs{A} = 1$ and for a set $R \su \S$ and $x \in \R$, we write $R - x \coeq \set{(y - x, u) \co (y, u) \in R}$.
Note that for translation-invariant~$f$, we can replace $f (\De_m - c(\De_m), \PP - c(\De_m), \PP' - c(\De_m))$ by $f(\De_m, \PP, \PP')$.
Note also that if~$f$ is bounded, \autoref{prop:finite_intensity} shows that the expectation~\eqref{eq:palm} is well-defined.

%
%
Now we can derive an integral representation of the Palm distribution.
We write $\bs{p}_m \coeq (p_1, \dots, p_m)$ for an $m$-tuple of points, and for $u \in [0, 1]$, we introduce the notations $\bs{p}_m(u) \coeq ((0, u), p_1, \dots, p_m)$ and $\vpp_m(u) \coeq \set{(0, u), p_1, \dots, p_m}$ for the tuple and corresponding set, respectively.

\begin{proposition}[Distribution of the typical $m$-simplex]\label{prop:palm_distribution}
    Let $m \ge 0$, $\g < 1$, $\g' < 1/(m + 1)$, and let $f \co (\R \ti [0, 1])^{m + 1} \ti \Nl \ti \Nl \to \R_+$ be an arbitrary nonnegative measurable functional depending on an $m$-simplex as well as on the point processes~$\PP$ and~$\PP'$.
    Then,
    \[ \E[\big]{f(\De_m^\ast, \PP, \PP')} = \f{\la^{m + 1}}{\la_m(m + 1)!} \int_{[0, 1] \ti \S^m} \E[\big]{f(\vpp_m(u), \PP \cup \vpp_m(u), \PP') \1{\vpp_m(u) \in \Sigma_m}} \d \pp_m \d u. \]
    Here, if $\vpp_m(u)$ does not consist of precisely $m + 1$ elements, then we let $f(\vpp_m(u), \PP \cup \vpp_m(u) , \PP') \coeq 0$.
\end{proposition}

As a first application of the Palm distribution, we show that the higher-order degree distributions are scale-free.
In classical binary networks, the degree of a vertex equals the number of edges it is incident to.
For instance, in the bipartite graph~$\Gb$, the degree of an author-vertex equals the number of papers this author has written.
Hence, similarly, for an $(m + 1)$-element subset of authors, we define its higher-order degree as the number of papers this subset has collaborated on.
With $\mrm{Fin}(\PP)$ denoting the family of finite subsets of~$\PP$, we define the bipartite graph with vertex sets $\mrm{Fin}(\PP)$ and~$\PP'$, and there is an edge between a simplex $\De \in \mrm{Fin}(\PP)$ and $p' \in \PP'$ if and only if $\max_{p \in \De} F(p, p'; \g, \g') \le \b$.
The degree $\deg(\De)$ of a simplex~$\De$ is then the number of points in~$\PP'$ to which all points in~$\De$ are connected to:
\[ \deg(\De) \coeq \PP'(B_F(\De, \b)) \coeq \#(\PP' \cap B_F(\De, \b)). \]
Note that $\deg(\De)$ depends in fact on~$\Gb$ and not only on~$\Gh$.
In the context of the Dowker complex, $\deg(\De)$ can also be interpreted as the number of witnesses for the $m$-simplex~$\De$.
We now show that the typical higher-order degrees are scale-free, characterized by a power-law distribution.
In particular, for $m = 0$, we recover the classical degree distribution of the authors.

%
%
\begin{theorem}[Scale-freeness of higher-order degrees]\label{thm:deg}
    Let $m \ge 0$, $\g < 1$, $\g' < 1/(m + 1)$.
    Then,
    \[ \lim_{k \tff} \f{\log[\big]{\P[\big]{\deg (\De_m^\ast) \ge k}}}{\log{k}} = m - \f{m + 1}{\g}. \]
\end{theorem}

We note that given a subset~$\s$ of~$\PP$ and a subset~$\t$ of~$\PP'$, one can say that $(\s, \t)$ forms a biclique in the bipartite graph~$\Gb$ if every $\PP$-vertex in~$\s$ is connected to every $\PP'$-vertex in~$\t$ and vice versa.
Then, the degree of~$\s$ is the number of bicliques of the form $(\s, \t)$ where~$\t$ is a set of cardinality exactly one.

From a topological perspective, we note that \autoref{thm:deg} is loosely related to the multicover bifiltration.
Here, considering a union of balls, a point is $k$-covered if it is contained in at least~$k$ of the balls.
Similarly, in \autoref{thm:deg} we have $\deg(\De_m\ph{\ast}) \ge k$ if there are at least~$k$ points in the set $B_F(\De_m, \b)$.

One of the key advantages of a stochastic network model is that, when working with data, we can statistically test whether the model is a good fit for the considered data.
To carry out this approach rigorously, we aim to use test statistics that reflect the key topological properties of the Dowker complex associated with~$\Gh$.

This opens the door towards considering invariants from topological data analysis, such as the Betti numbers.
Thus, we start by establishing the asymptotic normality of the Betti numbers in the regime $\g < 1/4$ in the window $\S_n \coeq [0, n] \ti [0, 1]$.
To make this precise, fix $m \ge 0$ and let $\b\ph{(n)}_m$ denote the $m$th Betti number of $\Gh(\PP \cap \S_n, \PP' \cap \S_n)$.
Moreover, let $\NN(0, \s^2)$ denote the normal distribution with mean~$0$ and variance~$\s^2$.

%
%
\begin{theorem}[Asymptotic normality of Betti numbers]\label{thm:bet}
    Let $m \ge 0$.
    Let $\g < 1/4$ and $\g' < 1/(4(m + 1))$.
    Then, in distribution,
    \[ n^{-1/2} \bigl( \b_m^{(n)} - \E[\big]{\b_m^{(n)}} \bigr) \xrightarrow[n \tff]{d} \NN(0, \s^2) \qquad \text{for some $\s^2 \ge 0$}. \]
\end{theorem}
Note that instead of letting $n \to \ff$, one could alternatively let~$\la$,~$\la'$, and~$\b$ tend to infinity.

Note also that by the Dowker duality theorem~\citep{dowker}, $\b_m\ph{(n)}$ coincides with the Betti number that is obtained from the dual Dowker complex, where the roles of~$\PP$ and~$\PP'$ are reversed.

While the asymptotic normality of test statistics is highly convenient in applied statistics, we stress that \autoref{thm:bet} requires that $\g < 1/4$.
Taking into account \autoref{thm:deg} for $m = 0$, we see that this means that the typical degree distribution has a finite fourth moment.
However, when considering complex networks, we often encounter the situation where even the variance is infinite.
Hence, it is no longer reasonable to expect a normal distribution in the limit, since the latter exhibits light tails.
In analogy to the classical setting of sums of independent heavy-tailed random variables, we expect that the suitably recentered and rescaled distribution converges to a stable distribution.

However, since the Betti numbers are a highly refined topological quantity, giving a rigorous proof of the stable limit convergence in the delicate heavy-tailed setting is difficult.
We also note that while there is an ample variety of asymptotically normal test statistics for spatial network models~\citep{shirai,yukCLT}, the stable setting has been considered so far only in selected isolated cases~\citep{dst}.

While establishing a stable limit for the Betti numbers seems to be out of reach for the moment, similarly to \citet{hj23}, we can prove a stable limit for a much simpler test statistic, namely the edge count
\[ S_n \coeq \sum_{P_i \in \PP \cap \S_n} \deg(P_i), \]
where we recall that $\deg(P_i)$ denotes the degree of a point~$P_i$ in the bipartite graph~$\Gb$.

%
%
\begin{theorem}[Normal and stable limits of edge counts]\label{thm:stab}
    Let $\g' < 1/3$.
    Then, the following distributional limits hold as $n \to \ff$.
    \begin{enumerate}[label=(\alph*)]
        \item
            Let $\g \in (0, 1/2)$.
            Then, $n^{-1/2}(S_n - \E{S_n}) \mathrel{\xrightarrow[n \tff]{d}} \NN(0, \s^2)$ for some $\s^2 > 0$.
            \label{thm:stab:thin}
        \item
            Let $\g \in (1/2, 1)$.
            Then, $n^{-\g}(S_n - \E{S_n}) \mathrel{\xrightarrow[n \tff]{d}} \SS_{\g^{-1}}$, where $\SS_{\g^{-1}}$ is a $\g^{-1}$-stable random variable.
            \label{thm:stab:heavy}
    \end{enumerate}
\end{theorem}

While \autoref{thm:stab} establishes rigorously the desired normal and stable limits in the asserted regimes, it is difficult to apply to hypothesis tests in an actual data analysis.
This is because the parameter~$\b$ is typically tuned such that the expected number of edges in the bipartite model matches the quantity observed in the data.

Hence, to have practically useful test statistics, we now define
\[ S_{n, m} \coeq \# \set[\big]{\De_m \in \Sigma_m \co c(\De_m) \in [0, n]} \]
as the number of $m$-simplices centered in $[0, n]$.

The fact that we can extend the proof of \autoref{thm:stab} to the case of $m$-simplex count is remarkable.
Indeed, in the simpler model considered in \citet{hj23}, this functional seemed out of reach.
The explanation is that for the task of proving central and stable limit theorems, our bipartite model is surprisingly more accessible than the one from \citet{hj23}.
The reason is that in the present model, the crucial covariance between simplex counts in disjoint regions becomes accessible since we can condition on the set $\PP'$ and then apply the formula for total covariance.
These essential variance and covariance computations are summarized in the following auxiliary result.

%
%
\begin{theorem}[Normal and stable limits of simplex counts]\label{thm:simp}
    Let $\g' < 1 / (2m + 1)$.
    Then, the following distributional limits hold as $n \to \ff$.
    \begin{enumerate}[label=(\alph*)]
        \item
            Let $\g \in (0, 1/2)$.
            Then, $n^{-1/2}(S_{n, m} - \E{S_{n, m}}) \mathrel{\xrightarrow[n \tff]{d}} \NN(0, \s^2)$ for some $\s^2 > 0$.
            \label{thm:simp:thin}
        \item
            Let $\g \in (1/2, 1)$.
            Then, $n^{-\g}(S_{n, m} - \E{S_{n, m}}) \mathrel{\xrightarrow[n \tff]{d}} \SS_{\g^{-1}}$, where $S_{\g^{-1}}$ is a $\g^{-1}$-stable random variable.
            \label{thm:simp:heavy}
    \end{enumerate}
\end{theorem}

To relate our model to existing concepts in the literature, we note that a standard construction exists to retrieve a hypergraph from a bipartite graph.
This could also be done in our setting with the bipartite graph~$\Gb$.
After taking closures under subsets, we recover the hypergraph~$\Gh$ introduced above.

In the remainder of the paper, the parameter $\b$ is held fixed.
Hence, to ease notation, in the rest of the paper, we write $B_F(\vpp_m)$ instead of $B_F(\vpp_m, \b)$.

\section[Palm distribution]{Proof of Propositions~\ref{prop:finite_intensity} and~\ref{prop:palm_distribution}}\label{sec:palm_distribution}

%
%
\begin{proof}[Proof of \autoref{prop:finite_intensity}]
    As before, let $A \su \R$ be an arbitrary Borel set with $\abs{A} = 1$.
    Furthermore, let $g_m(p_0, \dots, p_m, \PP')$ denote the indicator of the event that an $m$-tuple of points forms an $m$-simplex and that the marks of $p_0, \dots,p_m$ are ordered ascending.
    The $m$-simplex intensity $\la_m$ is then given by
    \[ \begin{aligned}
        \la_m &= \E[\Big]{\sum_{\De_m \in \Sigma_m} \1[\big]{c(\De_m) \in A}} = \la^{m + 1} \int_{(A \ti [0, 1]) \ti \S^m} \E[\big]{g_m \bigl( (p_0, p_1, \dots, p_m), \PP' \bigr)} \d(p_0, \pp_m) \\
        &= \la^{m + 1} \int_{[0, 1] \ti \S^m} \E[\big]{g_m(\pp_m(u), \PP')} \d \pp_m \d u,
    \end{aligned} \]
    where we used the Mecke formula \citep[Theorem~4.4]{poisBook} in the first step, and integrated with respect to the location of~$p_0$ in the second step.
    Note that the number of points in the common neighborhood $B(\vpp_m(u))$ of the vertices is Poisson-distributed.
    Hence, by the Markov inequality, the probability that this set is nonempty is given by:
    \[ \la_m = \f{\la^{m + 1}}{(m + 1)!} \int_0^1 \int_{\S^m} \P[\big]{\PP'(B(\vpp_m(u))) \ge 1} \d \pp_m \d u \le \f{\la^{m + 1} \la'}{(m + 1)!} \int_0^1 \int_{\S^m} \abs{B(\vpp_m(u))} \d \pp_m \d u. \]
    The inner integral is given by
    \[ \int_{\S^m} \abs{B(\vpp_m(u))} \d \pp_m = \int_\S \1[\big]{\abs{z} \le \b u^{-\g} w^{-\g'}} \int_{\S^m} \1[\big]{\abs{z - y_i} \le \b v_i^{-\g} w^{-\g'} \co 1 \le i \le m} \d \pp_m \d p'. \]
    Next, we integrate with respect to the location coordinates $\bs{y} \coeq (y_1, \dots, y_m)$, the marks $\bs{v} \coeq (v_1, \dots, v_m)$ as well as $z, w, u$.
    This gives for the above
    \[ \begin{aligned}
        &\int_\S \1[\big]{\abs{z} \le \b u^{-\g} w^{-\g'}} \int_{[0, 1]^m} \prod_{i = 1}^m \bigl( 2 \b v_i^{-\g} w^{-\g'} \bigr) \d \bs{v} \d p' = \Bigl( \f{2 \b}{1 - \g} \Bigr)^m \int_\S \1[\big]{\abs{z} \le \b u^{-\g} w^{-\g'}} w^{-m \g'} \d p' \\
        &\qquad = \f{(2 \b)^{m+1}}{(1 - \g)^m} u^{-\g} \int_0^1 w^{-(m+1) \g'} \d w = \f{(2 \b)^{m + 1}}{(1 - \g)^m} \f{u^{-\g}}{1 - (m + 1) \g'},
    \end{aligned} \]
    where we have used that $\g < 1$ and $\g' < 1 / (m+1)$.
    Using again that $\g < 1$, we conclude that $\la_m < \ff$.
\end{proof}

%
%
\begin{proof}[Proof of \autoref{prop:palm_distribution}]
    Let $A \su \R$ be an arbitrary Borel set with $\abs{A} = 1$.
    Let the lowest-mark $\PP$-vertex $c(\De_m)$ of a set of $m + 1$ points $\De_m$ be denoted by $p_0 \coeq (x, u)$.
    The expectation of a nonnegative functional~$f$ is given by~\eqref{eq:palm}.
    By another application of the Mecke formula, and by integrating over~$x$, we obtain
    \[ \E[\big]{f(\De_m^\ast, \PP, \PP')} = \f{\la^{m+1}}{\la_m} \int_{[0, 1] \ti \S^m} \E[\Big]{f \bigl( \vpp_m(u), \PP \cup \vpp_m(u), \PP' \bigr) g_m(\pp_m(u), \PP')} \d \pp_m \d u, \]
    as asserted.
\end{proof}

\section[Scale-freeness of higher-order degrees]{Proof of \autoref{thm:deg}}\label{sec:deg}

Since the proofs of the upper and lower bounds are very different, we deal with them in Sections~\ref{ssec:deg_upper} and~\ref{ssec:deg_lower}, separately.
We start with two lemmas that are used in both the proof of the upper and lower bounds.

\begin{lemma}[Pairwise intersections]\label{lem:pairwise_intersection}
    Let $0 \le u \le v \le 1$, and set $o \coeq (0, u) \in \S$ and $p \coeq (y, v) \in \S$.
    Then, the measure of the intersection of their neighborhood $\abs{B(\set{o, p})}$ can be upper bounded as follows:
    \[ \abs{B(\set{o, p})} \le \f{2 \b}{1 - \g'} v^{- \g} s_\w (u, y)^{1 - \g'} \qquad \text{where} \qquad s_\w (u, y) \coeq \bigl( 2 \b u^{-\g} \abs{y}^{-1} \bigr)^{{1}/{\g'}} \w 1. \]
\end{lemma}
\begin{proof}
    For a point $(z, w) \in \S$ to be connected to both $(0, u)$ and $(y, v)$, it must hold that $\abs{z} \le \b u^{-\g} w^{-\g'}$ and that $\abs{y - z} \le \b v^{-\g} w^{-\g'}$.
    As $u \le v$, we find from the triangle inequality that $\abs{y} \le 2 \b u^{-\g} w^{-\g'}$ and hence $w \le (2 \b u^{-\g} \abs{y}^{-1})^{1/\g'}$.
    On the other hand, $w \le 1$, so $w \le s_\w(u, v)$ and using only the second condition above for~$z$,
    \[ \abs{B(\set{o, p})} \le \int_0^{s_\w(u, y)} \int_{-\b v^{-\g} w^{-\g'}}^{\b v^{-\g} w^{-\g'}} \d z \d w = \int_0^{s_\w(u, y)} 2 \b v^{-\g} w^{-\g'} \d w = \f{2 \b}{1 - \g'} v^{-\g} s_\w(u, y)^{1 - \g'}, \]
    as asserted.
\end{proof}

\begin{lemma}[Upper bound of $\int F(s_\w(u, y)) \d y$]\label{lem:int_f_s_dy}
    Let $\r > \g'$ and $F \co (0, \ff) \to [0, 1]$ be an integrable function.
    Assume that $c \coeq \sup_{x > 0} x^{-\r}F(x) < \ff$.
    Then, for all $a \ge 0$,
    \[ \int_{-\ff}^\ff F(s_\w(u, y)) \d y \le 2 a + \f{2 c \g'}{\r - \g'} (2 \b)^{\r / \g'} a^{1 - \r/\g'} u^{- \r \g / \g'}. \]
\end{lemma}
\begin{proof}
    We use the upper bound of~$F$ and the symmetry of~$\abs{y}$ to write
    \[ \begin{aligned}
        \int_{-\ff}^\ff F(s_\w(u, y)) \d y &\le 2 \int_0^a 1 \d y + 2 c \int_a^\ff s_\w(u, y)^\r \d y \le 2 a + 2 c \bigl( 2 \b u^{-\g} \bigr)^{\r / \g'} \int_a^\ff y^{-\r / \g'} \d y \\
        &\le 2 a + \f{2 c \g'}{\r - \g'} (2 \b)^{\r / \g'} a^{1 - \r / \g'} u^{- \r \g / \g'},
    \end{aligned} \]
    where we have used that $\r > \g'$.
\end{proof}

%
%
\subsection{Proof of \autoref{thm:deg}---upper bound}\label{ssec:deg_upper}

The key idea for the proof of the upper bound is to note that if all points of $\De_m$ have a common neighbor in~$\PP'$, then also all pairs of $\PP$-vertices in $\De_m$ have a common neighbor.
This observation will be used to subsequently bound the integrals in the Palm probabilities.

%
%
\begin{proof}[Proof of the upper bound.]
    By \autoref{prop:palm_distribution} we obtain that
    \[ \P{\deg(\De_m^\ast) \ge k} = \f{\la^{m+1}}{\la_m(m+1)!} \int_{[0, 1]} \int_{\S^m} \P{\PP'(B(\vpp_m(u))) \ge k} \d \pp_m \d u. \]
    Let $b \coeq \la' \abs{B(\vpp_m(u))}$ and note that $X \coeq \PP'(B(\vpp_m(u)))$ is Poisson-distributed with parameter~$b$.
    Next, we bound the probability by distinguishing whether $b \ge k/2$:
    \[ \P{X \ge k} \le \1{b \ge k/2} + \P{X \ge k} \1{b < k/2}. \]
    Our goal is to show that the second term is negligible compared to the first for large values of~$k$.
    First,
    \[ \P{X \ge k} = \P{X \ge 1} \f{\P{X \ge k}}{\P{X \ge 1}}. \]
    Let $b_\vee \coeq 2 \vee e^2 b$.
    The numerator is upper bounded by $(b / k)^{k/2}$ whenever $k \ge b_\vee$ using \citet[Lemma~1.2]{geomgraphs}.
    If $b < 1/2$, the denominator can be bounded from below by $1 - \exp{- b} \ge b / 2$.
    Otherwise, the denominator can be bounded from below by a finite constant $1/c_2$.
    Thus, for $k \ge b_\vee$, we have
    \[ \f{\P{X \ge k}}{\P{X \ge 1}} \le c_1 \Bigl( \f{b}{k} \Bigr)^{k/2} \Bigl( \f{2}{b} \1[\Big]{b < \f{1}{2}} + c_2 \1[\Big]{b \ge \f{1}{2}} \Bigr) \le c_3 \Bigl( \f{b}{k} \Bigr)^{k/2 - 1}
    \]
    for some $c_1, c_2, c_3 > 0$.
    Considering the indicator $\1{b < k/2}$ in the integrand, we can conclude that
    \[ \int_{[0, 1] \ti \S^m} \P{X \ge k} \1[\big]{b < k / 2} \d (u, \pp_m) \le c_3 2^{1 - k / 2} \int_{[0, 1] \ti \S^m} \P{X \ge 1} \d (u, \pp_m) = c_4 2^{-k / 2}, \]
    where $c_4 > 0$ and in the last step we used that the integral above is finite due to \autoref{prop:finite_intensity}.
    As shown below, the $\1{b \ge k / 2}$ term decays as a power law with increasing~$k$.
    Thus,
    \[ \limsup_{k \tff} \f{\log[\big]{\P{\deg(\De_m^\ast) \ge k}}}{\log{k}} \le \limsup_{k \tff} \f{\displaystyle \log[\Big]{\int_{[0, 1] \ti \S^m} \1[\big]{b \ge k/2} \d(u, \pp_m)}}{\log{k}}. \]
    In particular, for $m = 0$,
    \[ \begin{aligned}
        \int_{[0, 1]} \1[\big]{b \ge k/2} \d u &= \int_{[0, 1]} \1[\Big]{\la' \int_{\R \ti [0, 1]} \1[\big]{\abs{z} \le \b u^{-\g} w^{-\g'}} \d (z, w) \ge k/2} \d u \\
        &= \int_{[0, 1]} \1[\bigg]{\f{2 \b \la'}{1 - \g'} u^{-\g} \ge \f{k}{2}} \d u \in O( k^{-1/\g}).
    \end{aligned} \]
    This concludes the proof for $m = 0$.

    From now on, we assume $m \ge 1$.
    Our goal is to upper bound the expression
    \[ \limsup_{k \tff} \f{\log[\big]{\P{\deg( \De_m^\ast) \ge k}}}{\log{k}} \le \limsup_{k \tff} \Bigl( \f{1}{\log{k}} \log[\bigg]{\int_{[0, 1] \ti \S^m} \1[\Big]{b \ge \tf{k}{2}} \d (u, \pp)} \Bigr). \]
    We assume that the points~$p_i$ are ordered in increasing order of their marks, i.e., $u \le v_1 \le \cdots \le v_m$.
    Note that
    \[ b \le \la' \abs{B(\vpp_1(u))} \w \min_{i = 1, \dots, m - 1} \la' \abs{B(\set{p_i, p_{i + 1}})}. \]
    Due to translation invariance and \autoref{lem:pairwise_intersection} we find that for some constant $c > 0$,
    \[ \la' \abs{B(\set{p_i, p_{i + 1}})} = \la' \abs{B(\set{(0, v_i), (y_{i + 1} - y_i, v_{i + 1})})} \le c v_{i + 1}^{-\g} s_\w (v_i, y_{i + 1} - y_i)^{1 - \g'}. \]
    Using this upper bound,
    \[ \int_{[0, 1] \ti \S^m} \1[\big]{b \ge k / 2} \d(u, \pp_m) \le \int_{[0, 1] \ti \S^m} \prod_{i = 0}^{m - 1} \1[\Big]{c v_{i + 1}^{-\g} s_\w (v_i, y_{i + 1} - y_i)^{1 - \g'} \ge k} \d(u, \pp_m), \]
    where we have set $v_0 \coeq u$ and $y_0 \coeq 0$.
    In the indicators, we can express the upper limit of the marks~$v_i$,
    \[ \int_0^{c k^{- 1 / \g}} \int_{\S^m} \prod_{i = 0}^{m - 1} \1[\Big]{v_{i + 1} \le c s_\w (v_i, y_{i + 1} - y_i)^{(1 - \g') / \g} k^{-1/\g}} \d \pp_m \d u, \]
    where we set the upper limit of the integral with respect to~$u$ to $c k^{- 1 / \g}$ as $u \le v_1 \le c k^{- 1 / \g}$.
    We now substitute $y_{i + 1}' \coeq y_{i + 1} - y_i$ for $i = 0, \dots, m - 1$.
    Furthermore, let $p_i' \coeq (y_i', v_i)$ denote the $i$th point with the new coordinates.
    Then, we have
    \[ \int_0^{c k^{- 1 / \g}} \int_{\S^m} \prod_{i = 0}^{m - 1} \1[\Big]{v_{i + 1} \le c s_\w(v_i, y_{i + 1}')^{(1 - \g')/\g} k^{-1/\g}} \d \pp_m' \d u. \]
    Note that the point $p'_m = (y'_m, v_m)$ only appears in one of the indicators in the product.
    We integrate this indicator with respect to~$p'_m$.
    Then, for large~$k$,
    \[ \int_{\R \ti [0, 1]} \1[\Big]{v_m \le c s_\w(v_{m - 1}, y_m')^{(1 - \g') / \g} k^{-1/\g}} \d p'_m = \int_\R c k^{-1 / \g} s_\w(v_{m - 1}, y'_m)^{(1 - \g') / \g} \d y'_m \in O(k^{-1/\g} v_{m - 1}^{-\g}), \]
    where in the first step we noted that the indicator represents an upper bound of~$v_m$ and integrated it out, and in the second step we used \autoref{lem:int_f_s_dy} with $\r = (1 - \g') / \g$ and $a = u^{-\g}$.
    Note that the finiteness of the integral requires that $\g < (1 - \g') / \g'$.
    We can follow the same steps for integrating out with respect to $p'_{m - 1}$.
    The only difference is that now a $v_{m - 1}^{-\g}$ appears in the integration, resulting in a $k\ph{(1 - \g) / \g}$ factor.
    Then,
    \[ \begin{aligned}
        \int_{[0, 1] \ti \S^m} \1[\big]{b \ge k/2} \d(u, \pp_m) &\le c k^{- {(1 + (m - 2)(1 - \g))}/\g} \int_0^{c k^{- 1 / \g}} \int_{\R \ti [0, 1]} \1[\Big]{ c v_1^{-\g} s_\w(u, y'_1)^{1 - \g'} \ge k} v_1^{-\g} \d p'_1 \d u \\
        &\le c k^{- {(1 + (m - 1)(1 - \g))}/\g} \int_0^{c k^{- 1 / \g}} u^{-\g} \d u \in O(k^{m - {(m + 1)}/\g}).
    \end{aligned} \]
    This leads to
    \[ \lim_{k \tff} \biggl( \f{1}{\log{k}} \log[\bigg]{\int_{[0, 1] \ti \S^m} \1[\big]{b \ge k/2} \d(u, \pp_m)} \biggr) = m - \f{m + 1}{\g}, \]
    as asserted.
\end{proof}

%
%
\subsection{Proof of \autoref{thm:deg}---lower bound}\label{ssec:deg_lower}

The idea for proving the lower bound is to construct specific configurations that lead to higher-order degrees exceeding~$k$.
We then show that these configurations occur with a sufficiently high probability.

%
%
\begin{proof}[Proof of the lower bound.]
    Let $R_k \coeq [0, k] \ti [0, (\b/k) ^{1/\g}] \su \S$ and note that for $p \coeq (y, v) \in R_k$ and $p' \coeq (z, w) \in [0, k] \ti [0, 1]$,~$p$ and~$p'$ are connected.
    Indeed, $\abs{y - z} \le k \le \b u^{-\g} w^{-\g'}$.
    This gives for $\nu(k) \coeq \la' k - \log{2}$ by \autoref{prop:palm_distribution}
    \[ \begin{aligned}
        \P{\deg(\De_m^\ast) \ge \nu(k)} &\ge \f{\la^{m + 1}}{\la_m (m + 1)!} \int_{R_k^m} \int_{[ 0, (\b/k)^{1/\g}]} \P[\big]{\deg(\vpp_m(u)) \ge \nu(k)} \d u \d \bs{p}_m \\
        &\ge \f{\la^{m + 1}(\b / k) ^{1/\g}}{\la_m (m + 1)!} \int_{R_k^m} \P[\Big]{\PP'([0, k] \ti [0, 1]) \ge \nu(k)} \d \bs{p}_m.
    \end{aligned} \]
    As $\PP' ([0, k] \ti [0, 1]) \sim \Poi{\la' k}$, its median can be lower bounded by $\nu(k)$ \citep[Theorem~2]{poimedian}.
    Thus, lower bounding the integrand, we obtain
    \[ \P{\deg(\De_m^\ast) \ge \nu(k)} \ge \f{\la^{m+1}(\b/k) ^{1/\g}}{2 \la_m (m+1)!} \abs{R_k}^m = \f{\la^{m + 1} \b^{(m + 1) / \g}}{2 \la_m(m + 1)!} k^{m - (m + 1)/\g}. \]
    Hence,
    \[ \liminf_{k \tff} \f{\log[\big]{\P{\deg( \De_m^\ast) \ge k}}}{\log{k}} = \liminf_{k \tff} \f{\log[\big]{\P{\deg( \De_m^\ast) \ge \nu(k)}}}{\log{k}} \ge m - \f{m + 1}{\g} \]
    as asserted.
\end{proof}

\section{Proof of \autoref{thm:bet}}\label{sec:bet}

The proof of the CLT for the Betti numbers relies on the general Poisson CLT by \citet[Theorem~3.1]{yukCLT}.
This result requires verification of two conditions: a stabilization condition and a moment condition.

The stabilization condition follows from an adaptation of the arguments provided in \citet{hj23} and \citet{shirai}.
However, the moment condition is more involved than in \citet{hj23} or in \citet{shirai}.
This is because, in contrast to \citet{hj23} and \citet{shirai}, the existence of an $m$-simplex can no longer be determined just from the knowledge of the positions of the $m + 1$ $\PP$-vertices.
Indeed, the existence fundamentally involves the second Poisson process~$\PP'$.

Let $\b(\PP, \PP') \coeq \b_{n, m}(\PP, \PP')$ be the $m$th Betti number of $\Gh_n(\PP \cap \S_n, \PP' \cap \S_n)$.
We introduce the two special points $o \coeq (0, U), o' \coeq (0, W)$.
Furthermore, let $\Gh_{n, o} \coeq \Gh ((\PP \cup \set{o}) \cap \S_n, \PP' \cap \S_n)$ denote the hypergraph in the window~$\S_n$ with the additional $\PP$-vertex~$o$, and similarly, let $\Gh_{n, o'}(\PP \cap \S_n, (\PP' \cup \set{o'}) \cap \S_n)$ denote the hypergraph in the window~$\S_n$ with the additional $\PP'$-vertex~$o'$.
Using these notations, the add-one cost operators
\[ \begin{aligned}
    \de(\Gh_n, o ) &\coeq \b(\Gh_{n, o }) - \b(\Gh_n) \\
    \de(\Gh_n, o') &\coeq \b(\Gh_{n, o'}) - \b(\Gh_n).
\end{aligned} \]

To apply \citet[Theorem~3.1]{yukCLT}, we need to verify the following conditions:
\begin{itemize}
    \item \textbf{moment condition:} $\sup_{n \ge 1} \E{\de(\Gh_n, o)^4} < \ff$ and \\
    $\sup_{n \ge 1} \E{\de(\Gh_n, o')^4} < \ff$
    \item \textbf{weak stabilization:} $\lim_{n \tff} \de(\Gh_n, o) < \ff$ and $\lim_{n \tff} \de(\Gh_n, o') < \ff$
\end{itemize}

\begin{proof}[Proof of the moment condition]
    By \citet[Lemma~2.9]{shirai} \citet[Chapter~VII]{edHar}, $\de(\Gh_n, o)$ and $\de(\Gh_n, o')$ are bounded above by the number of $m$- and $(m + 1)$-simplices containing the additional $\PP$-vertex $o = (0, U) \in \PP$ or $o' = (0, W) \in \PP'$.

    We first consider the case of $\de(\Gh_n, o')$.
    The number of new $m-1$-simplices is given as the number of $m$-tuples of points in the neighborhood $B(o')$ of the new point $(0, W) \in \PP'$:
    \[ \E[\big]{\de(\Gh_n, o')^4} \le \E[\bigg]{\biggl( \binom{\PP(B(o'))}{m} \biggr)^4} \le \E[\big]{\PP(B(o'))^{4 m}}. \]
    The number of points in the neighborhood $\PP(B(o'))$ is Poisson distributed with mean $\la \abs{B(o')} = (2 \b \la) / (1 - \g) w^{-\g'}$, as shown in \autoref{prop:finite_intensity}.
    Since the $r$th factorial moment of a Poisson random variable with parameter $\mu > 0$ equals $\mu^r$, the $r$th moment itself is bounded above by $c_{\msf{M}, r} \mu^r$ for some constant $c_{\msf{M}, r} > 0$.
    Noting that $\g' < 1 / (4m)$, we conclude that
    \[ \E[\big]{\PP (B(o'))^{4m}} \le c_{\msf{M}, 4m} \E[\bigg]{\Bigl( \f{2 \b}{1 - \g} W^{-\g'} \Bigr)^{4m}} \le c \int_0^1 w^{- 4 m \g'} \d w < \ff \]
    for some $c > 0$.

    Next, we show that $\E{\de(\Gh_n, o)^4} < \ff$.
    The number of new $m$-simplices are formed by the sets of~$m$ points $\vPp_m \coeq \set{P_1, \dots, P_m} \in \PP^m_{\ne}$ for which a common neighbor $P' \in \PP'$ exists with the new $\PP$-vertex~$o$.
    As this number upper bounds $\de(\Gh_n, o)$, we have
    \[ \E{\de(\Gh_n, o)^4} \le \E[\bigg]{\biggl( \sum_{\bs{P}_m \in \PP_{\ne}^m} \1[\Big]{\PP' \bigl( B \bigl( \set{o} \cup \vPp_m \bigr) \bigr) \ge 1} \biggr)^4} \le \E[\bigg]{\biggl( \sum_{P' \in \PP'} \biggl( \binom{\PP(B(P'))}{m} \biggr) \1[\big]{P' \in B(o)} \biggr)^4}, \]
    where $\vPs_m \coeq \set{P_1', \dots, P_m'}$, and in the second step each term of the sum represents the number of $m$-tuples of $\PP$-vertices in the neighborhood of a point $P' \in \PP'$, and the indicator ensures that this point is in turn connected to the new $\PP$-vertex~$o$.
    We bound the binomial coefficient by $\PP(B(P'))^m$, and expand the fourth power of the sum by using the multinomial theorem to obtain
    \begin{alignat*}{2}
        \E{\de(\Gh_n, o)^4} \le {} &c_1 \E[\Big]{\sum_{\substack{\bs{P}'_4 \in \PP_{\ne}'^4}} \prod_{i \le 4} \PP \bigl( B(P'_i) \bigr)^m \1[\big]{\vPs_4 \su B(o)}} && \rule[0.5ex]{0.3cm}{0.2pt} \quad \text{term I} \\
            &+ c_2 \E[\Big]{\sum_{\bs{P}'_3 \in \PP_{\ne}'^3} \PP \bigl( B(P'_1) \bigr)^{2m} \prod_{i = 2}^3 \PP \bigl( B(P'_i) \bigr)^m \1[\big]{\vPs_3 \su B(o)}} && \rule[0.5ex]{0.3cm}{0.2pt} \quad \text{term II} \\
            &+ c_3 \E[\Big]{\sum_{\bs{P}'_2 \in \PP_{\ne}'^2} \PP \bigl( B(P'_1) \bigr)^{2m} \PP \bigl( B(P'_2) \bigr)^{2m} \1[\big]{\vPs_2 \su B(o)}} && \rule[0.5ex]{0.3cm}{0.2pt} \quad \text{term III} \\
            &+ c_4 \E[\Big]{\sum_{\bs{P}'_2 \in \PP_{\ne}'^2} \PP \bigl( B(P'_1) \bigr)^{3m} \PP \bigl( B(P'_2) \bigr)^m \1[\big]{\vPs_2 \su B(o)}}    && \rule[0.5ex]{0.3cm}{0.2pt} \quad \text{term IV} \\
            &+ c_5 \E[\Big]{\sum_{\substack{P' \in \PP'}} \PP \bigl( B(P') \bigr)^{4m} \1[\big]{P' \in B(o)}} && \rule[0.5ex]{0.3cm}{0.2pt} \quad \text{term V}
    \end{alignat*}
    for some $c_1, \dots, c_5 > 0$.
    To show the moment bound, it is enough to show that each of the above terms is finite.
    We begin with using the multivariate Mecke formula for term I:
    \[ I \coeq c_1 \E[\Big]{\sum_{\substack{\bs{P}'_4 \in \PP_{\ne}'^4}} \prod_{i \le 4} {\PP \bigl( B(P'_i) \bigr)^m} \1[\big]{\vPs_4 \su B(o)}} = c_1 \iint_{[0, 1] \ti \S^4} \E[\Big]{\prod_{i \le 4} \PP \bigl( B(p'_i) \bigr)^m \1[\big]{p_i' \in B(o)}} \d u \d \pp_4'. \]
    Introducing the notation $p'_i \coeq (z_i, w_i)$ for the coordinates of the points~$p'_i$, the indicator represents an upper bound for the coordinates $\abs{z_i} \le \b u^{-\g} w_i^{-\g'}$.
    Next, we use H\"older's inequality to upper bound the expectation of the product:
    \[ I \le 2^4 c_1 \iint_{[0, 1] \ti [0, 1]^4} \prod_{i \le 4} \biggl( \int_0^{\b u^{-\g} w_i^{-\g'}} \E[\Big]{\PP \Bigl( B((z_i, w_i)) \Bigr)^{4m}}^{1/4} \d z_i \biggr) \d u \d \bs{w}. \]
    Hence,
    \[ I \le 2^4 c_1 \iint_{[0, 1] \ti [0, 1]^4} \prod_{i \le 4} \Bigl( \int_0^{\b u^{-\g} w_i^{-\g'}} c_{\msf{M}, 4m} \Bigl( \f{2 \b}{1 - \g} w_i^{-\g'} \Bigr)^m \d z_i \Bigr) \d u \d \bs{w}. \]
    Note that the upper bound does not depend on the variable~$z_i$, as it is translation invariant.
    Also, we can upper bound the exponential term by a sufficiently large constant~$C$ to get
    \[ I \le \Bigl( \f{C (2 \b)^{(m + 1)}}{(1 - \g)^m} \Bigr)^4 c_1 \int_0^1 u^{- 4 \g} \d u \int_{[0, 1]^4} \prod_{i = 1}^4 w_i^{- (m + 1) \g'} \d \bs{w}. \]
    The integral over the~$w_i$ is finite whenever $\g' < 1/(m+1)$, and the integral over~$u$ is finite if $\g < 1/4$.
    The other terms can be bounded similarly, and we obtain a finite bound for $\E{\de(\Gh_n, o)^4}$.
    The conditions for the finiteness of the integrals are given as follows:
    \begin{itemize}
        \item term I:   $\g < 1/4$ and $\g' < 1 / (  m + 1)$,
        \item term II:  $\g < 1/3$ and $\g' < 1 / (2 m + 1)$,
        \item term III: $\g < 1/2$ and $\g' < 1 / (2 m + 1)$,
        \item term IV:  $\g < 1/2$ and $\g' < 1 / (3 m + 1)$,
        \item term V:   $\g < 1$   and $\g' < 1 / (4 m + 1)$.
    \end{itemize}
    Combining these conditions, we obtain that $\g < 1/4$ and $\g' < 1 / (4 m + 1)$ are sufficient conditions for the finiteness of $\E{\de(\Gh_n, o)^4}$.
\end{proof}

%
%
\begin{proof}[Proof of the weak stabilization]
    In the following, we extend the arguments from \citet[Proposition~5.4]{shirai}.

    For $n \ge 1$ we write
    \[ \b_n \coeq \dim(Z_n) - \dim(B_n) \coeq \dim(Z(\Gh_n)) - \dim(B(\Gh_n)) \]
    for the Betti number of $\Gh_n$, where~$Z_n$ is the corresponding cycle space and~$B_n$ is the boundary space.

    %
    %
    We set for $o \coeq (0, U)$,
    \[ \b_{n, o} \coeq \dim(Z_{n, o}) - \dim(B_{n, o}) \coeq \dim(Z(\Gh_{n, o})) - \dim(B(\Gh_{n, o})). \]
    Hence, it suffices to prove the weak stabilization with respect to $\dim(Z_n)$ and $\dim(B_n)$ separately.
    Since the arguments are very similar, we henceforth only deal with the case $\dim(Z_n)$.
    We do this by showing that $\dim(Z_{n, o}) - \dim(Z_n)$ is non-decreasing and bounded in~$n$.

    For the boundedness, note that $\dim(Z_{n, o}) - \dim(Z_n) \le \deg_{m, n}(o)$, where $\deg_{m, n}(o)$ is the number of $m$-simplices in~$\S_n$ containing the typical $\PP$-vertex~$o$.
    This is because the $m$-simplices constructed from $\Gh_{n, o}$ can be decomposed into the set of $m$-simplices containing~$o$ and into the family of all $m$-simplices formed in $\Gh_n$ \citep[see][Lemma~2.9]{shirai}.
    To show that $\dim(Z_{n, o}) - \dim(Z_n)$ is non-decreasing, take $n_1 \le n_2$, and consider the canonical map
    \[ Z_{n_1, o} \to Z_{n_2, o} / Z_{n_2} \]
    with kernel $Z_{n_1, o} \cap Z_{n_2} \su Z_{n_1}$.
    We claim that the kernel is equal to $Z_{n_1}$.
    Let~$M_n$ and $M_{n, o}$ denote the set of $m$-simplices in $\Gh_n$ and $\Gh_{n, o}$, respectively, and consider an $m$-simplex $\s \in M_{n_1, o} \cap M_{n_2}$ forming an $m$-cycle~$z$.
    If $\s \in M_{n_1}$, then both $\s \in M_{n_1, o}$ and $\s \in M_{n_2}$.
    Therefore, $Z_{n_1} \su Z_{n_1, o'} \cap Z_{n_2}$, and the induced map
    \[ Z_{n_1, o} / Z_{n_1} \to Z_{n_2, o} / Z_{n_2} \]
    is injective.
    In particular, $\dim(Z_{n_1, o}) - \dim(Z_{n_1}) \le \dim(Z_{n_2, o}) - \dim(Z_{n_2})$, which shows the assertion.

    %
    %
    Now, we show the weak stabilization with respect to the typical $\PP'$-vertex~$o' \in \PP'$.
    As before, we write
    \[ \b_n' \coeq \dim(Z_{n, o'}) - \dim(B_{n, o'}) \coeq \dim(Z(\Gh_{n, o'})) - \dim(B(\Gh_{n, o'})), \]
    and prove weak stabilization with respect to $\dim(Z_n)$.

    Let~$M_n$ and $M_{n, o'}$ denote the set of $m$-simplices in $\Gh_n$ and $\Gh_{n, o'}$, respectively.
    Following the same arguments as above, the difference $\dim(Z_{n, o'}) - \dim(Z_n)$ is bounded from above by the number of $m$-simplices $\abs{M_{n, o'}} - \abs{M_n}$ appearing in $\Gh_n$ due to the addition of the typical $\PP'$-vertex~$o'$:
    \[ \dim(Z_{n, o'}) - \dim(Z_n) \le \abs{M_{n, o'}} - \abs{M_n}. \]
    In turn, $\abs{M_{n, o'}} - \abs{M_n}$ is bounded from above by the number of $m$-tuples of points $p \in \PP$ in the neighborhood $\PP(B(o'))$ of $o'$:
    \[ \dim(Z_{n, o'}) - \dim(Z_n) \le \biggl( \binom{\PP(B(o'))}{m} \biggr) < \ff, \]
    since the neighborhood $\PP(B(o'))$ is almost surely finite.

    We consider again the canonical map
    \[ Z_{n_1, o'} \to Z_{n_2, o'} / Z_{n_2}, \]
    where $n_1 \le n_2$.
    As $Z_{n_1} \su Z_{n_2}$, the kernel of the map is $Z_{n_1, o'} \cap Z_{n_2}$.
    If $Z_{n_1} = Z_{n_1, o'} \cap Z_{n_2}$, we can conclude the proof as above.

    In contrast to the case above, however, $\dim(Z_{n, o'}) - \dim(Z_n)$ is not necessarily monotone, since $Z_{n_1, o'} \cap Z_{n_2} \ne Z_{n_1}$, i.e., there can be cycles in the kernel that are not in $Z_{n_1}$.
    To see this, consider a cycle $z \in Z_{n_1, o'} \sm Z_{n_1}$.
    If~$z$ contains a simplex formed by a $\PP'$-vertex in the increased window size~$n_2$, but not in the window size~$n_1$, then $z \in Z_{n_2} \sm Z_{n_1}$.
    Thus, $z \notin Z_{n_1}$, but is in the kernel of the map $Z_{n_1, o'} \to Z_{n_2, o'} / Z_{n_2}$.
    To solve this problem, we construct a random~$N_2$ such that for any two window sizes $n_2 \ge n_1 \ge N_2$ we have $Z_{n_1, o'} \cap Z_{n_2} = Z_{n_1}$ in the kernel of the map $Z_{n_1, o'} \to Z_{n_2, o'} / Z_{n_2}$ the sequence $(\dim(Z_{n, o'}) - \dim(Z_n))$ is monotone.
    Once we have shown this, we can conclude the proof as above.

    To construct~$N_2$, first observe that as $\PP(B(o'))$ is almost surely finite, there a random $N_1 \ge 1$ such that for all $n \ge N_1$, the number of $m$-simplices $\abs{M_{n, o'}} - \abs{M_n}$ is a non-increasing sequence.
    As all $\PP$-vertices in the neighborhood $B(o')$ have finite neighborhoods almost surely, there also exists $N_2 \ge N_1 \in \R$ such that for all $n \ge N_2$, the number of $m$-simplices $\abs{M_{n, o'}} - \abs{M_n}$ is a non-decreasing sequence, and thus constant for such~$n$.
    Then, we conclude the proof as above.
\end{proof}

\section[Normal and stable limits of the edge count]{Proof of \autoref{thm:stab}}

The proof idea for \autoref{thm:stab} is similar to the edge-count CLTs and stable limits obtained in \citet{hj23}.
For the CLT, we use a general result on associated random variables, whereas for the stable case, we use a truncation argument and a comparison with the independent case.
Besides some moment computations of the degree of a typical node, the main difficulty in the proof is to approximate the random number of nodes in the interval $[0, n]$ by a deterministic number.

%
%
We begin by proving an auxiliary result that will be used in both the normal and the stable cases.
Recall the definition of the neighborhood $B(\De)$ of a set of points~$\De$ as it was defined in~\eqref{eq:ball}.
As before, for $u \le 1$ and $\vpp_m \su \PP$, we write $B(\vpp_m(u))$ for the common neighborhood of $\set{(0, u), p_1, \dots, p_m}$.

%
%
\begin{lemma}[Scaling of $B(\vpp_m(u))$]\label{lem:mu2m}
    Let $u \le 1$, $m \ge 1$, $0 < \g < 1$ and $0 < \g' < 1/(m+1)$.
    It holds that
    \[ \int_{\S^m} \abs{B(\vpp_m(u))} \d \pp_m = \f{u^{-\g}}{1-(m + 1) \g'} \f{(2 \b)^{m + 1}}{(1 - \g)^m}. \]
\end{lemma}
\begin{proof}
    We have that
    \[ \begin{aligned}
        \int_{\S^m} \abs{B(\vpp_m(u))} \d \pp_m &= \iint_{\S \ti \S^m} \1[\big]{p' \in B(\vpp_m(u))} \d \pp_m \d p' \\
        &= \int_\S \1[\big]{(z, w) \in B(o)} \biggl( \int_\S \1[\big]{\abs{z - y} \le \b v^{-\g} w^{-\g'}} \d(y, v) \biggr)^m \d(z, w) \\
        &= \int_\S \1[\big]{(z, w) \in B(o)} \Bigl( \f{2 \b w^{-\g'}}{1 - \g} \Bigr)^m \d(z, w) \\
        &= u^{-\g} \int_0^1 \f{\bigl( 2 \b w^{-\g'} \bigr)^{m + 1}}{(1 - \g)^m} \d w = \f{u^{-\g}}{1 - (m + 1) \g'} \f{(2 \b)^{m + 1}}{(1 - \g)^m},
    \end{aligned} \]
    where we used the notations $p \coeq (y, v)$ and $p' \coeq (z, w)$.
\end{proof}

%
%
\subsection{Proof of \autoref{thm:stab}---normal limit}\label{sec:clt}

The idea of the proof is to apply the CLT \citet[Theorem~4.4.3]{whitt} which holds for stationary sequences of identically distributed associated random variables $\bs{T} \coeq T_1, T_2, \dots$, and requires that $\sum_{k \ge 1} \Cov{T_1, T_k} < \ff$.
Recall that the sequence of random variables $\bs{T}$ is associated if and only if $\Cov{f(T_1, \dots, T_k), g(T_1, \dots, T_k)} \ge 0$ for all non-decreasing functions~$f$,~$g$ for which $\E{f(T_1, \dots, T_k)}$, $\E{g(T_1, \dots, T_k)}$, $\E{f(T_1, \dots, T_k) g(T_1, \dots, T_k)}$ exist \citep[Definition~1.1]{associatedrvs}.
\begin{proof}[Proof of \autoref{thm:stab}~\ref{thm:stab:thin}]
    Let $\g < 1/2$ and define for $i \ge 1$,
    \[ T_i \coeq \sum_{P_j \in \PP \cap ([i - 1, i] \ti [0, 1])} \deg(P_j). \]
    Since the degree $\deg(x, u)$ is increasing in the Poisson process~$\PP'$, we conclude from the Harris-FKG inequality \citep[Theorem~20.4]{poisBook} that the sequence $T_1, \dots, T_k$ is associated.
    For $A \su \R$ let
    \[ S(A) \coeq \sum_{P_j \in \PP \cap A \ti [0, 1]} \deg(P_j). \]
    Then, the Mecke equation gives
    \[ \begin{aligned}
        \Var{S(A)} &= \E[\Big]{\sum_{P_i \ne P_j \in \PP \cap A} \deg(P_i) \deg(P_j)} + \E[\Big]{\sum_{P_i \in \PP \cap A} \deg(P_i)^2} - \E[\Big]{\sum_{P_i \in \PP \cap A}\deg(P_i)}^2 \\
        &= \la^2 \iint_{(A \ti [0, 1])^2} \E{\deg(p_1) \deg(p_2)} \d p_1 \d p_2 \\
        &\phantom{=} + \la \int_{A \ti [0, 1]} \E{\deg(p_1)^2} \d p_1 - \la^2 \biggl( \int_{A \ti [0, 1]} \E{\deg(p_1)} \d p_1 \biggr)^2 \\
        &= \la \abs{A} \int_0^1 \E{\deg(0, u)^2} \d u + \la^2 \iint_{(A \ti [0, 1])^2} \Cov{\deg(p_1), \deg(p_2)} \d p_1 \d p_2,
    \end{aligned} \]
    where the last equality holds due to translation invariance of $\deg(\any)$.

    To bound the first integral, note that for all $u \in [0, 1]$, the random variable $\deg(0, u)$ is Poisson distributed with parameter $(2 \b \la' u^{-\g}) / (1 - \g')$.
    Hence, $\E{\deg(0, u)^2} = \f{2 \b \la'}{1 - \g'} u^{-\g} (1 + \f{2 \b \la'}{1 - \g'} u^{-\g}) \in O(u^{-2 \g})$.
    Thus, we obtain for $\g < 1/2$ that
    \[ \int_0^1 \E{\deg(0, u)^2} \d u < \ff. \]

    Next, we deal with the second term and note that the covariance in the integrand is given by
    \[ \begin{aligned}
        \Cov{\deg(p_1), \deg(p_2)} &= \E[\Big]{\sum_{P' \in \PP'} \1[\big]{P' \in B(\set{p_1, p_2})}} + \E[\Big]{\sum_{P'_1 \ne P'_2 \in \PP'} \1[\big]{P'_1 \in B(p_1), P'_2 \in B(p_2)}} \\
        &\phantom{=} - \E[\Big]{\sum_{P' \in \PP'} \1[\big]{P' \in B(p_1)}} \E[\Big]{\sum_{P' \in \PP'} \1[\big]{P' \in B(p_2)}} \\
        &= \la' \int_\S \1[\big]{p' \in B(\set{p_1, p_2})} \d p' = \la' \abs{B(\set{p_1, p_2})}.
    \end{aligned} \]
    From \autoref{lem:mu2m} with $m = 1$ and translation invariance, we obtain that
    \[ \int_{A \ti [0, 1]} \int_\S \abs{B(\set{p_1, p_2})} \d p_1 \d p_2 = \f{1}{1 - 2 \g'} \f{(2 \b)^2}{1 - \g} \int_0^1 u^{-\g} \d u < \ff. \]
    Hence, $\Var{T_1} = \Var{S([0, 1])} < \ff$.
    Similarly, we obtain that
    \[ \begin{aligned}
        \sum_{k \ge 2} \Cov{T_1, T_k} &= \E[\Big]{\sum_{\substack{P_i \in \PP \cap [0, 1] \ti [0, 1], \\
        P_j \in \PP \cap [1, \ff) \ti [0, 1]}} \deg(P_i) \deg(P_j)} - \E[\Big]{\sum_{P_i \in \PP \cap [0, 1] \ti [0, 1]} \deg(P_i)} \E[\Big]{\sum_{P_j \in \PP \cap [1, \ff) \ti [0, 1]} \deg(P_j)} \\
        &\le \la^2 \int_{[0, 1]^2} \int_\S \Cov{\deg(p_1), \deg(p_2)} \d p_2 \d p_1 \le \f{\la^2 \la'}{1 - 2 \g'} \Bigl( \f{2 \b}{1 - \g} \Bigr)^2 \int_0^1 u^{-\g} \d u < \ff,
    \end{aligned} \]
    thereby showing the finiteness of $\sum_{k \ge 1} \Cov{T_1, T_k}$.
\end{proof}

\subsection{Proof of \autoref{thm:stab}---stable limit}\label{sec:stab}

The idea of the proof of \autoref{thm:stab}~\ref{thm:stab:heavy} is to apply \citet[Theorem~4.5.2]{whitt}, which says the following.
Let $X_i, i \ge 1$ be i.i.d.\@ nonnegative random variables with $\P{X_1 > x} \sim A x^{-\a}$ for some $\a \in (1, 2)$ and $A > 0$.
Then $n^{-1/\a} \sum_{i \le n} (X_i - \E{X_i})$ converges in distribution to an $\a$-stable random variable.

%
%
\begin{proof}[Proof of \autoref{thm:stab}~\ref{thm:stab:heavy}]
   Let $u_n \coeq n^{-b}$ for an arbitrary choice of $b \in (2/3, 1)$.
   We decompose~$S_n$ into
   \[ S_n = \sna + \snb \coeq \sum_{P \in \PP \cap \S_{n, \ge u_n}} \deg(P) + \sum_{P \in \PP \cap \S_{n, \le u_n}} \deg(P), \]
   where $\S_{n, \ge u} \coeq [0, n] \ti [u, 1]$ and $\S_{n, \le u} \coeq [0, n] \ti [0, u]$.

   First, we show that $n^{-\g}(\sna - \E{\sna})$ converges to~$0$ in probability as $n \to \ff$.
   By Chebychev's inequality, this follows once we have shown that $\Var{\sna} \in o(n^{2 \g})$.
   Note that
   \begin{equation} \Var{\sna} = \la \int_{\S_{n, \ge u_n}} \E{\deg(p)^2} \d p + \la^2 \iint_{\S_{n, \ge u_n}^2} \Cov{\deg(p_1), \deg(p_2)} \d p_1 \d p_2 \label{varstabge} \end{equation}
   Here, by translation invariance, the first term is bounded by $(n u_n^{1 - 2 \g}) / (2 \g - 1)$, which is in $o(n^{2 \g})$ since $1 - b(1 - 2 \g) < 2 \g$ for all $b < 1$ and $\g > 1/2$.
   For the second term, \autoref{lem:mu2m} gives the bound
   \[ 2 n \la^2 \int_{u_n}^1 \int_{\S_{\ge u}} \abs{B(\set{(0, u), (y, v)})} \d(y, v) \d u \in O(n), \]
   where $\S_{\ge u} \coeq \R \ti [u, 1]$.
   Thus, $n^{-\g} (\sna - \E{\sna}) \xrightarrow[n \tff]{\mbb{P}} 0$.

   For $u \le 1$, we set $\mu(u) \coeq \E{\deg(0, u)} = (2 \b u^{-\g}) / (1 - \g')$ and
   \[ \snc \coeq \sum_{(X, U) \in \PP \cap \S_{n, \le u_n}} \mu(U) \qquad \text{and} \qquad \snd \coeq \sum_{i \le \ceil{\la n}} \mu(U_i) \1[\big]{U_i \le u_n}, \]
   where $U_1, U_2, \dots$ are i.i.d.\@ uniforms on $[0, 1]$.
   We next prove that $\E{\abs{\snb - \snc}} + \E{\abs{\snc - \snd}} \in o(n^\g)$.
   For the first expectation, note that, by the Mecke equation,
   \[ \E{\abs{\snb - \snc}} \le n \la \int_0^{u_n} \E{\abs{\deg(0, u) - \mu(u)}} \d u \le n \la \int_0^{u_n} \Var{\deg(0, u)}^{1/2} \d u = n \la \int_0^{u_n} \mu(u)^{1/2} \d u, \]
   Therefore, $\E{\abs{\snb - \snc}} \in O(n u_n^{1 - \g/2})$.
   Since $1 - b (1 - \g/2) < \g$, which holds for $b > 2/3$, we conclude that $n^{-\g} \E{\abs{\snb - \snc}} \to 0$.

   Second, we obtain a Poisson random variable~$N$ with parameter $\la n$,
   \[ \E{\abs{\snc - \snd}} \le \E{\abs{N - \ceil{\la n}}} \int_0^{u_n} \mu(u) \d u \le \bigl( \Var{N}^{1/2} + 1 \bigr) \int_0^{u_n} \mu(u) \d u. \]
   Now, $\Var{N} = \la N$ and the integral above is in $O(u_n^{1 - \g})$.
   Therefore, $\E{\abs{\snc - \snd}} \in O(n^{1/2} u_n^{1 - \g})$.
   Since $1/2 - b (1 - \g) < \g$, which holds for all~$b$ if $\g > 1/2$, it follows that $n^{-\g} \E{\abs{\snc - \snd}} \to 0$.

   It therefore suffices to prove that $n^{-\g} (\snd - \E{\snd})$ converges in distribution to a stable random variable.
   Note that by \citet[Lemma~11]{hj23},
   \[ \sum_{i \le \ceil{\la n}} \mu(U_i) \1[\big]{U_i > u_n} \xrightarrow[n \tff]{L^2} 0. \]
   Finally, \citet[Theorem~4.5.2]{whitt} implies that $n^{-\g^{-1}} \sum_{i \le \ceil{\la n}} (\mu(U_i) - \E{\mu(U_i)})$ converges in distribution to a $\g^{-1}$-stable random variable.
\end{proof}

\section[Normal and stable limits of the simplex count]{Proof of \autoref{thm:simp}}\label{sec:stabm}

We now extend the proof of \autoref{thm:stab} to the case of $m$-simplices.
More precisely, we consider
\[ S_{n, m} \coeq \sum_{P \in \PP \cap \S_n} \d_m(P), \]
where for $p = (x, u)$, we set
\[ \d_m(p) \coeq \f{1}{m!} \sum_{(P_1, \dots, P_m) \in (\PP \cap \S_{> u})_{\ne}^m} \1[\big]{\PP' \cap B(\set{p, P_1, \dots, P_m}) \ne \es} \]
is the number of $m$-simplices containing~$p$.

%
%
\begin{lemma}[Moment computations]\label{lem:condcov}
   Let $m \ge 1$.
   Then,
   \begin{enumerate}[label=(\alph*)]
      \item $\E{\d_m(0, u)^2} \in O(u^{-2 \g})$. \label{lem:condcov:a}
      \item $\E{\Var{\d_m(0, u) \given \PP'}} \in O(u^{1 - 3 \g})$. \label{lem:condcov:b}
      \item $\int_\S \abs{\E{\Cov{\d_m(0, u), \d_m(p) \given \PP'}}} \d p \in O(u^{1 - 3 \g})$. \label{lem:condcov:c}
      \item $\Var{\E{\d_m(0, u) \given \PP'}} \in O(u^{-\g})$. \label{lem:condcov:d}
      \item $\int_\S \abs{\Cov{\E{\d_m(0, u) \given \PP'}, \E{\d_m(p) \given \PP'}}} \d p \in O(u^{-\g})$. \label{lem:condcov:e}
   \end{enumerate}
\end{lemma}

We defer the technical computations to the end of this section and first explain how to conclude the proof of \autoref{thm:simp} based on \autoref{lem:condcov}.
For $\g < 1/2$, we follow the same strategy as in the proof of \autoref{thm:stab}~\ref{thm:stab:thin}.

\begin{proof}[Proof of \autoref{thm:simp}~\ref{thm:simp:thin}]
   Let
   \[ T_i \coeq \sum_{P_j \in [i - 1, i] \ti [0, 1]} \d_m(P_j). \]
   We apply \citet[Theorem~4.4.3]{whitt} which holds for sequences of identically distributed associated random variables $T_1, T_2, \dots$, and need to show that $\sum_{k \ge 1} \Cov{T_1, T_k} < \ff$.

   By the formula for the total variance, we have $\Var{T_1} = \E{\Var{T_1 \given \PP'}} + \Var{\E{T_1 \given \PP'}}$, where by \autoref{lem:condcov}~\ref{lem:condcov:a},~\ref{lem:condcov:c}, and translation invariance
   \[ \E{\Var{T_1 \given \PP'}} = \la \int_{\S_1} \E{\d_m(p)^2} \d p + \la^2 \iint_{\S_1^2} \E{\Cov{\d_m(p_1), \d_m(p_2) \given \PP'}} \d p_1 \d p_2 < \ff. \]
   Moreover, by \autoref{lem:condcov}~\ref{lem:condcov:e} and translation invariance,
   \[ \Var{\E{T_1 \given \PP'}} = \la^2 \iint_{\S_1^2} \Cov{\E{\d_m(p_1) \given \PP'}, \E{\d_m(p_2) \given \PP'}} \d p_1 \d p_2 < \ff. \]

   Analogously, we obtain that
   \[ \begin{aligned}
      \sum_{k \ge 2} \Cov{T_1, T_k} &= \sum_{k \ge 2} \Bigl( \E{\Cov{T_1, T_k \given \PP'}} + \Cov[\big]{\E{T_1 \given \PP'}, \E{T_k \given \PP'}} \Bigr) \\
      &= \la^2 \int_{\S_1} \int_{[1, \ff) \ti [0, 1]} \E{\Cov{\d_m(p_1), \d_m(p_2) \given \PP'}} \d p_1 \d p_2 \\
      &\phantom{=} + \la^2 \int_{\S_1} \int_{[1, \ff) \ti [0, 1]} \Cov[\big]{\E{\d_m(p_1) \given \PP'}, \E{\d_m(p_2) \given \PP'}} \d p_1 \d p_2 < \ff.
   \end{aligned} \]
\end{proof}

Next, we prove \autoref{thm:simp}~\ref{thm:simp:heavy}.
\begin{proof}[Proof of \autoref{thm:simp}~\ref{thm:simp:heavy}]
   The idea is to proceed analogously to the proof of \autoref{thm:stab}~\ref{thm:stab:heavy}.
   We assume that $\g' < 1 / (2 m + 1)$, and choose $u_n = n^{- b}$, where $b \in (2/3, 1)$, similarly to the proof of \autoref{thm:stab}.
   Then, we split $S_{n, m}$ into
   \[ S_{n, m} = S_{n, m}^\ge + S_{n, m}^\le \coeq \sum_{P_i \in \PP \cap \S_{n, \ge u_n}} \d_m(P_i) + \sum_{P_i \in \PP \cap \S_{n, \le u_n}} \d_m(P_i). \]

   We now show that $n^{-\g}(S_{n, m}^\ge - \E{S_{n, m}^\ge})$ converges to~$0$ in probability.
   This follows, as soon as we have proved that $\Var{S_{n, m}^\ge} \in o(n^{2 \g})$.
   We have
   \[ \Var{S_{n, m}^\ge} = n \la \int_{u_n}^1 \E{\d_m(0, u)^2} \d u + n \la^2 \int_{\S_{1, \ge u_n}} \int_{\S_{n, \ge u_n}} \Cov{\d_m(x, u), \d_m(p)} \d p \d(x, u). \]
   By \autoref{lem:condcov}~\ref{lem:condcov:a},~\ref{lem:condcov:c}, and translation invariance, there are constants $c_1, c_2 > 0$ such that the above is bounded by
   \[ c_1 n \la \int_{u_n}^1 u^{-2 \g} \d u + c_2 n \la^2 \int_{u_n}^1 u^{1 - 3 \g} \d u \in O(n u_n^{1 - 2 \g}), \]
   which is in $o(n^{2 \g})$, since $1 - b (1 - 2 \g) < 2 \g$ for $\g > 1/2$.
   Thus, $n^{-\g}(S_{n, m}^\ge - \E{S_{n, m}^\ge}) \stackrel{\mbb{P}}{\longrightarrow} 0$.
   Next, we let
   \[ S_{n, m}^{(1)} \coeq \sum_{(X, U) \in \PP \cap \S_n} \mu_m(U), \qquad S_{n, m}^{(2)} \coeq \sum_{i \le \ceil{\la n}} \mu_m(U_i) \1[\big]{U_i \le u_n}, \]
   where $\mu_m(u) \coeq \E{\d_m(0, u)}$, and $U_1, U_2, \dots$ are i.i.d.\@ uniform random variables on $[0, 1]$.
   We next prove that $\E{\abs{S_{n, m}^\le - S_{n, m}\ph{(1)}}} + \E{\abs{S_{n, m}\ph{(1)} - S_{n, m}\ph{(2)}}} \in o(n^\g)$.
   For the first expectation, note that by the Mecke equation,
   \[ \E{\abs{S_{n, m}^\le - S_{n, m}^{(1)}}} \le n \la \int_0^{u_n} \E{\abs{\d_m(0, u) - \mu_m(u)}} \d u \le n \la \int_0^{u_n} \Var{\d_m(0, u)}^{1/2} \d u. \]
   From \autoref{lem:condcov}~\ref{lem:condcov:b},~\ref{lem:condcov:d}, we have that $\Var{\d_m(0, u)} \in O(u^{1 - 3 \g})$.
   Hence,
   \[ \E{\abs{S_{n, m}^\le - S_{n, m}^{(1)}}} \le n u_n^{3 (1 - \g) / 2} \in o(n^\g), \]
   since $1 - b (3 (1 - \g) / 2) < \g$ for $b > 2/3$.
   From here, we conclude the proof as in \autoref{thm:stab}~\ref{thm:stab:heavy}.
\end{proof}

Finally, we conclude by proving \autoref{lem:condcov}.
We note that the proof of \autoref{lem:condcov:b} is largely parallel to \autoref{lem:condcov:c}, and the proof of \autoref{lem:condcov:d} is mainly parallel to \autoref{lem:condcov:e}.
Nevertheless, to make the manuscript self-contained, we provide the details.
\begin{proof}[Proof of \autoref{lem:condcov}]
   Henceforth, we use the abbreviations
   \[ \begin{aligned}
      \pp_{i, m} &\coeq (p_{m - i + 1}, \dots, p_{2 m - i}) \\
      \vpp_{i, m} &\coeq \set{p_{m - i + 1}, \dots, p_{2 m - i}}, \\
      \vpp_{i, m}(u) &\coeq \set{(0, u), p_{m - i + 1}, \dots, p_{2 m - i}}.
   \end{aligned} \]
   Furthermore, to ease notation, we introduce for a set~$A$ the abbreviation $\PP'_A \coeq \PP' \cap B(A)$.

   \paragraph{\autoref{lem:condcov:a}.}
   Conditioned on~$\PP'$, $\d_m(0, u)$ is a $U$-statistic.
   Therefore, by \citet[Lemma~3.5]{reitzner} with $k_i \coeq 1 / (i! ((m - i)!)^2)$,
   \[ \E{\d_m(0, u)^2 \given \PP'} = \sum_{i = 0}^m \la^{2 m - i} k_i \int_{\S_{\ge u}^{2 m - i}} \P[\big]{\PP'_{\vpp_m(u)} \ne \es, \PP'_{\vpp_{i, m}(u)} \ne \es \given \PP'} \d \pp_{2 m - i}. \]
   Therefore,
   \[ \begin{aligned}
      \E{\d_m(0, u)^2} &= \sum_{i = 0}^m \la^{2 m - i} k_i \int_{\S_{\ge u}^{2 m - i}} \P[\big]{\PP'_{\vpp_m(u)} \ne \es, \PP'_{\vpp_{i, m}(u)} \ne \es} \d \pp_{2 m - i} \\
      &\le \sum_{i = 0}^m \la^{2 m - i} k_i \int_{\S_{\ge u}^{2 m - i}} \P{\PP'_{\vpp_{2 m - i}(u)} \ne \es} + \P{\PP'_{\vpp_m(u)} \ne \es} \P{\PP'_{\vpp_{i, m}(u)} \ne \es} \d \pp_{2 m - i},
   \end{aligned} \]
   where we used that for all Borel sets $A, B \su \S$, $\P{\PP' \cap A \ne \es, \PP' \cap B \ne \es} = \P{\PP' \cap A \cap B \ne \es} + \P{\PP' \cap A \sm B \ne \es, \PP' \cap B \sm A \ne \es}$.
   Next, we use that for all Borel set $A \su \R \ti [0, 1]$, $\P{\PP' \cap A = \es} \le \la' \abs{A}$ and obtain that the above is bounded by
   \[ \sum_{k = m}^{2 m} \la^k \int_{\S_{\ge u}^k} \Bigl( \la' \abs{B(\vpp_k(u))} + (\la')^2 \abs{B(\vpp_m(u))} \abs{B(\set{(0, u), p_{m + 1}, \dots, p_k})} \Bigr) \d \pp_k, \]
   which is of order $O(u^{-2 \g})$ by \autoref{lem:mu2m}.

   \paragraph{\autoref{lem:condcov:b}.}
   By \citet[Lemma~3.5]{reitzner},
   \[ \Var{\d_m(0, u) \given \PP'} = \sum_{i = 1}^m \la^{2 m - i} k_i \int_{\S_{\ge u}^{2 m - i}} \1[\big]{\PP'_{\vpp_m(u)} \ne \es} \1[\big]{\PP'_{\vpp_{i, m}(u)} \ne \es} \d \pp_{2 m - i}. \]
   Therefore,
   \[ \E{\Var{\d_m(0, u) \given \PP'}} = \sum_{i = 1}^m \la^{2 m - i} k_i \int_{\S_{\ge u}^{2 m - i}} \P{\PP'_{\vpp_m(u)} \ne \es, \PP'_{\vpp_{i, m}(u)} \ne \es} \d \pp_{2 m - i}. \]
   Similarly, as in \autoref{lem:condcov:a}, we find that the above is bounded by
   \[ \begin{aligned}
      &\int_{\S_{\ge u}^{2 m - i}} \P{\PP'_{\vpp_{2 m - i}(u)} \ne \es} + \P{\PP'_{\vpp_m(u)} \ne \es} \P{\PP'_{\set{(0, u), p_{m - i + 1}, \dots, p_{2 m - i}}} \ne \es} \d \pp_{2 m - i} \\
      &\qquad \le \int_{\S_{\ge u}^{2 m - i}} \la' \abs{B(\vpp_{2 m - i}(u))} + (\la')^2 \abs{B(\vpp_m(u))} \abs{B(\set{p_m, \dots, p_{2 m - i}})} \d \pp_{2 m - i},
   \end{aligned} \]
   where the last inequality follows by the same argument as in \autoref{lem:condcov:a}.
   Recall that $\vpp_{i, m}(u) = \set{(0, u), p_{m - i + 1}, \dots, p_{2 m - i}}$.
   Next, we apply \autoref{lem:mu2m} to $\abs{B(\vpp_{2 m - 1}(u))}$, integrated with respect to $\pp_{2 m - 1}$, and to $\abs{B(\vpp_{1, m})}$, integrated with respect to $p_{m + 1}, \dots, p_{2 m - i}$.
   This gives
   \[ \la^m (\la')^2 \biggl( \sum_{i = 1}^m k_i \f{(2 \b)^{m - i + 1}}{(1 - \g)^{m-i}} \f{1}{1 - (m - i + 1) \g'} \biggr) \int_{\S_{\ge u}^m} v_m^{-\g} \abs{B(\vpp_m(u))} \d \pp_m + O(u^{-\g}), \]
   where $p_m \coeq (y_m, v_m)$.
   Moreover, for all $u \in [0, 1]$,
   \begin{equation} \begin{aligned}
      &\int_{\S_{\ge u}^m} v_m^{-\g} \abs{B(\vpp_m(u))} \d \pp_m \nonumber \\
      &\qquad = \int_\S \int_{\S_{\ge u}^m} v_m^{-\g} \1[\big]{\abs{z} \le \b w^{-\g'} u^{-\g}}  \prod_{i = 1}^m \1[\big]{\abs{z - y_i} \le \b w^{-\g'} v_i^{-\g}} \d((y_1, v_1), \dots, (y_m, v_m)) \d(z, w) \nonumber \\
      &\qquad \le (2 \b)^m \int_\S \int_{[0, 1]^{m - 1} \ti [u, 1]} w^{-m \g'} v_1^{-\g} \cdots v_{m - 1}^{-\g} v_m^{-2 \g} \1[\big]{\abs{z} \le \b w^{-\g'} u^{-\g}} \d(v_1, \dots, v_m) \d(z, w) \nonumber \\
      &\qquad = \f{(2 \b)^{m + 1} (u^{1 - 2 \g} - 1) u^{-\g}}{(1 - \g)^{m - 1}(2 \g - 1)} \int_0^1 w^{-(m + 1) \g'} \d w = \f{(2 \b)^{m + 1} (u^{1 - 2 \g} - 1) u^{-\g}}{(1 - \g)^{m - 1}(2 \g - 1)(1 - (m + 1) \g')},
   \end{aligned} \label{bouvmgam} \end{equation}
   which shows that $\E{\Var{\d_m(0, u) \given \PP'}} \in O(u^{1 - 3 \g})$.

   \paragraph{\autoref{lem:condcov:c}.}
   From the polarization identity $\Cov{X, Y} = 1/4 (\Var{X + Y} - \Var{X - Y})$ and \citet[Lemma~3.5]{reitzner} (applied to the $U$-statistics $\d_m(0, u) + \d_m(p)$ and $\d_m(0, u) - \d_m(p)$) we obtain that
   \[ \begin{aligned}
          &\int_\S \Cov{\d_m(0, u), \d_m(p) \given \PP'} \d p \\
          &\qquad = \sum_{i = 1}^m \la^{2 m - i} k_i \int_\S \int_{\S_{\ge u}^m \ti \S_{\ge v}^{m - i}} \hspace{-.5cm} \1[\big]{\PP'_{\vpp_m(u)} \ne \es} \1[\big]{\PP'_{\set{p, p_{m - i + 1}, \dots, p_{2 m - i}}} \ne \es} \d \pp_{2 m - i} \d p,
      \end{aligned} \]
   where $p \coeq (y, v)$.
   Therefore, analogously to \autoref{lem:condcov:b},
   \[ \begin{aligned}
       &\int_\S \abs[\Big]{\E[\big]{\Cov{\d_m(0, u), \d_m(p) \biggiven \PP'}}} \d p \\
       &\qquad \le \sum_{i = 1}^m \la^{2 m - i} k_i \int_\S \int_{\S_{\ge u}^m \S_{\ge v}^{m - i}} \la' \abs[\big]{B(\vpp_{2 m - i}(u))} + (\la')^2 \abs[\big]{B(\vpp_m(u))} \abs[\big]{B(\set{p, p_m, \dots, p_{2 m - i}})} \d \pp_{2 m - i} \d p.
   \end{aligned} \]
   Now, we apply \autoref{lem:mu2m} to $\abs{B(\vpp_{2 m - 1}(u))}$ (integrated with respect to $\pp_{2 m - 1}$) and to $\abs{B(\set{p, p_m, \dots, p_{2 m - 1}})}$ (integrated with respect to $(p, p_{m + 1}, \dots, p_{2 m - i})$).
   This gives
   \[ \la^m (\la')^2 \Bigl( \sum_{i = 1}^m \f{(2 \b)^{m - i + 2}}{(1 - \g)^{m - i + 1}} \f{k_i}{1 - (m - i + 2) \g'} \Bigr) \int_{\S_{\ge u}^m} v_m^{-\g} \abs{B(\vpp_m(u))} \d \pp_m + O(u^{-\g}). \]
   From here, $\int_\S \E{\Cov{\d_m(0, u), \d_m(p) \given \PP'}} \d p \in O(u^{1-3\g})$ follows analogously to~\eqref{bouvmgam} as asserted.

   \paragraph{\autoref{lem:condcov:d}.}
   We have that
   \[ \E{\d_m(0, u) \given \PP'} = \f{\la^m}{m!} \int_{\S_{\ge u}^m} \1[\big]{\PP'_{\vpp_m(u)} \ne \es} \d \pp_m, \]
   and therefore,
   \begin{equation}
      \Var{\E{\d_m(0, u) \given \PP'}} = \f{\la^{2 m}}{(m!)^2} \int_{\S_{\ge u}^{2 m}} \Cov[\big]{\1[\big]{\PP'_{\vpp_m(u)} = \es}, \1[\big]{\PP'_{\vpp_{0, m}(u)} = \es}} \d \pp_{2 m}.
      \label{varcov}
   \end{equation}
   Now, we use that $\Cov{XY, XZ} = \E{Y} \E{Z} \Var{X}$ for independent $X, Y, Z$.
   We let
   \[ \begin{aligned}
      X &\coeq \1[\big]{\PP'_{\vpp_{2 m}(u)} = \es}, \\
      Y &\coeq \1[\big]{\PP' \cap (B(\vpp_m(u)) \sm B(\vpp_{0, m}(u))) = \es}, \\
      Z &\coeq \1[\big]{\PP' \cap (B(\vpp_{0, m}(u)) \sm B(\vpp_m(u))) = \es},
   \end{aligned} \]
   and deduce that
   \begin{equation} \begin{aligned}
      \Var{\E{\d_m(0, u) \given \PP'}} &= \f{\la^{2 m}}{m!} \int_{\S_{\ge u}^{2 m}} \P[\big]{\PP' \cap (B(\vpp_m(u)) \cup B(\vpp_{0, m}(u))) = \es} \P[\big]{\PP'_{\vpp_{2 m}(u)} \ne \es} \d \pp_{2 m} \nonumber \\
      &\le \f{\la^{2 m} \la'}{m!} \int_{\S_{\ge u}^{2 m}} \abs{B(\vpp_{2 m}(u))} \d \pp_{2 m}, \label{eq:msimcovar}
   \end{aligned} \end{equation}
   which shows by \autoref{lem:mu2m} that $\Var{\E{\d_m(0, u) \given \PP'}} \in O(u^{-\g})$.

   \paragraph{\autoref{lem:condcov:e}.}
   Analogously to~\eqref{varcov}, we find that
   \[ \begin{aligned}
       &\int_\S \abs[\big]{\Cov{\E{\d_m(0, u) \given \PP'}, \E{\d_m(p) \given \PP'}}} \d p \\
       &\qquad = \la^{2 m} \int_\S \int_{\S_{\ge u}^m \ti \S_{\ge v}^m} \abs[\Big]{\Cov[\big]{\1[\big]{\PP'_{\vpp_m(u)} \ne \es}, \1[\big]{\PP'_{\set{p} \cup \vpp_{i, m}} \ne \es}}} \d \pp_{2 m} \d p,
   \end{aligned} \]
   where $p = (y, v)$.
   Here, we apply the relation $\Cov{XY, XZ} = \E{Y} \E{Z} \Var{X}$ with
   \[ \begin{aligned}
      X &\coeq \1[\big]{\PP'_{\set{p} \cup \vpp_{2 m}(u)} = \es}, \\
      Y &\coeq \1[\big]{\PP' \cap (B(\vpp_m(u)) \sm B(\set{p} \cup \vpp_{0, m})) = \es}, \\
      Z &\coeq \1[\big]{\PP' \cap (B(\set{p} \cup \vpp_{0, m}) \sm B(\vpp_m(u))) = \es},
   \end{aligned} \]
   and deduce analogously to~\eqref{eq:msimcovar} that
   \[ \int_\S \abs[\big]{\Cov{\E{\d_m(0, u) \given \PP'}, \E{\d_m(p) \given \PP'}}} \d p \le \la^{2 m} \int_\S \int_{\S_{\ge u}^m \ti \S_{\ge v}^m} \abs{B(\set{p} \cup \vpp_{2 m}(u))} \d \pp_{2 m} \d p, \]
   which is by \autoref{lem:mu2m} of order $O(u^{-\g})$.
\end{proof}

\section{Simulation study}\label{sec:sim}

To demonstrate how the theoretical results can be applied to networks of finite size, this section presents a simulation study analyzing finite networks and Palm distributions.

To analyze the statistical properties of the model, we utilize a Monte Carlo approach: we implemented a \Cpp algorithm to generate several higher-order networks with identical model parameters.
The simulation of a finite network consists of the following steps:
\begin{enumerate}
    \item
        We first construct a torus of dimension~$1$, and generate two Poisson point processes~$\PP$,~$\PP'$.
        For practical reasons, we do not keep the intensity measures~$\la$ and~$\la'$ of the Poisson point processes constant throughout the simulations, but instead generate the network samples on a torus of Lebesgue measure~$1$, and adjust the Poisson intensity measures accordingly to generate networks with different expected number of $\PP$- and $\PP'$-vertices.
        After generating~$\PP$ and~$\PP'$, the vertices are equipped with a mark uniformly distributed in $[0, 1]$.
    \item
        Then, we generate the connections between the $\PP$- and $\PP'$-vertices using the connection rule described in~\eqref{eq:wrchm}.
        By using periodic boundary conditions, we can treat the vertices equivalently, irrespective of their position in the torus.
        Note that the parameter~$\b$ can be tuned to adjust the expected number of edges in the network.
        The expected $\De_0$-degree of a point $o = (0, u)$ is given by
        \[ \E[\Big]{\sum_{p' \in \PP'} \1[\big]{p' \in B(o)}} = 2 \b \la' / ({(1 - \g) (1 - \g')}), \]
        where we followed the same arguments as in the proof of~\ref{lem:mu2m}, with $m = 0$.
    \item
        Finally, we construct the higher-order network on~$\PP$, where the set of $m$-simplices is determined according to~\eqref{eq:simplices}.
\end{enumerate}

\begin{figure} \centering
    \includegraphics{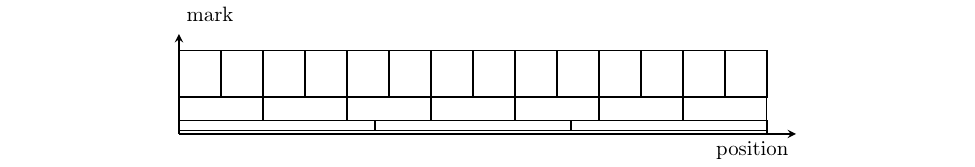}
    \caption{Partition of~$\S$ into rectangles.}\label{fig:rectangles}
\end{figure}

The simulation of the network is computationally expensive.
If we checked the connection condition for every~$\PP$- and $\PP'$-vertex pair, our algorithm would have a time complexity of $O(\la \la')$, where~$\la$ and~$\la'$ are the intensity of the Poisson processes for the $\PP$- and $\PP'$-vertices, respectively.
When deciding which points in~$\PP$ and~$\PP'$ are connected, we partition the space~$\S$ into rectangles containing multiple points, similarly to the quadtree data structures~\citep{quadtree}.
This way, we can reduce the computational complexity of the simulation as
\begin{itemize}
    \item if the furthest points of two rectangles with the highest marks connect, it is certain that all point pairs in the rectangles are connected;
    \item if the closest points of two rectangles with the lowest marks do not connect, it is certain that no point pairs in the rectangles are connected.
\end{itemize}
The width of the rectangles is set to be proportional to $b^{-\g}$ and $b^{-\g'}$ for the~$\PP$- and $\PP'$-vertices, respectively, where~$b$ denotes the bottom mark of the rectangle.
This results in a layout of rectangles as shown in \autoref{fig:rectangles}.

Before generating a finite network sample, we set the model parameters to be of the same order of magnitude when fitted to the analyzed datasets.
The intensity measures of the Poisson point processes are set to $\la = \la' = \num{100000}$.
Aligned with these parameters, the number of~$\PP$- and $\PP'$-vertices were \num{100335} and \num{93919}, respectively.
The~$\g$,~$\g'$ model parameters were set to $\g = 0.7$, $\g' = 0.2$.
Furthermore,~$\b$ was set so that the expected number of~$\De_0$-degrees in the infinite network limit was~$3$: $\b = 3.6e-5$.

First, we examine three degree distributions of the networks, which are analyzed in \autoref{thm:deg}.
\begin{itemize}
    \item{The $\De_0$-degree distribution characterizes the distribution of the number of $\PP'$-vertices connected to a typical $\PP$-vertex.}
    \item{The $\De_1$-degree distribution is the distribution of the number of $\PP'$-vertices connected to a pair of $\PP$-vertices.}
    \item{Finally, the distribution of the number of $\PP$-vertices connected to a typical $\PP'$-vertex is described by the $\De_0'$-degree distribution.}
\end{itemize}

\begin{figure} \centering
    \includegraphics{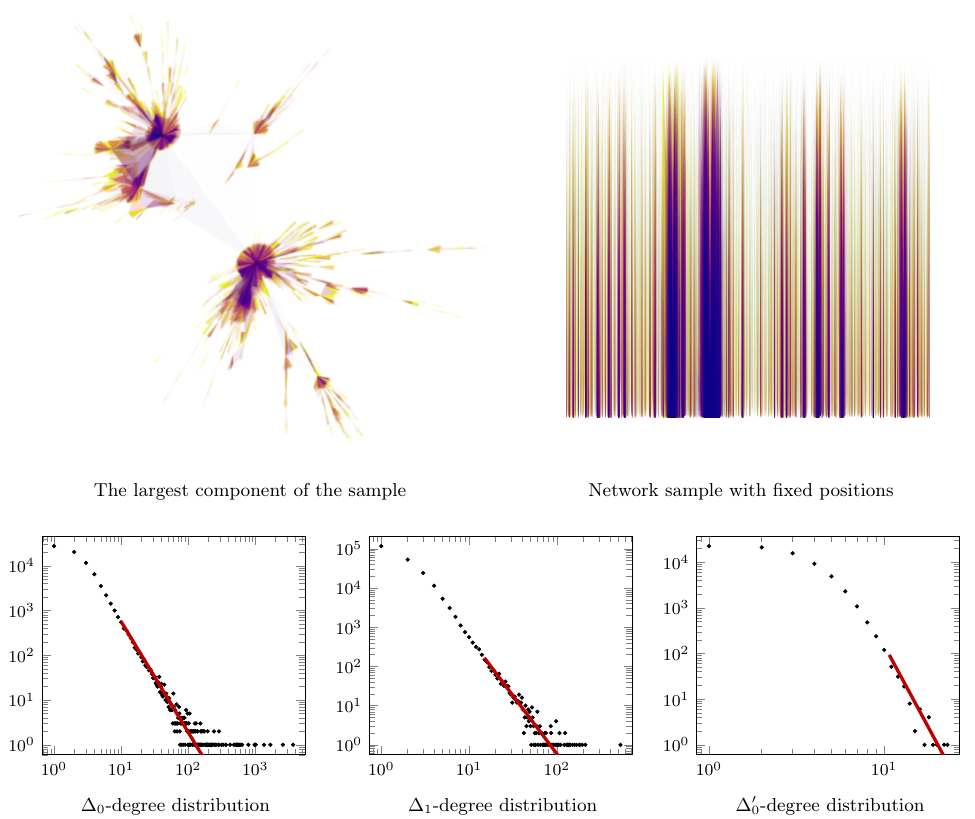}
    \caption{Main properties of a network sample generated by our model}\label{fig:model_sample}
\end{figure}

The main characteristics of a network generated with these parameters are shown in \autoref{fig:model_sample}.
In the top row of \autoref{fig:model_sample}, one can see the largest component of the network (left).
A polygon is drawn representing the convex hull of a set of $\PP$-vertices whenever they form a simplex in the network.
The higher number of $\PP$-vertices forms a simplex; the darker the color of the corresponding polygon.
From the largest component of the network, which contained \num{12385} $\PP$-vertices, one can see that the network is sparse, and that the network layout is dominated by a few dense regions corresponding to $\PP$-vertices with high degrees.

When looking at the positions of the $\PP$-vertices embedded in the space~$\S$ (right), one can see that the $\PP$-vertices form clusters that are connected by the above-mentioned heavy vertices.
As illustrated in the bottom row of~\ref{fig:model_sample}, the various degree distributions of the~$\PP$- and $\PP'$-vertices are all heavy-tailed, and the degree distribution of the $\PP$-vertices is heavier than the degree distribution of the $\PP'$-vertices, aligned with our theoretical results.

We also rely on Palm calculus to simulate typical simplices in networks devoid of finite-size effects.
To do so, we first place a typical $\PP$-vertex $(0, u)$ at the origin.
Then, we simulate $\PP'$-vertices $(y_1, v_1), \dots, (y_n, v_n) \in B(0, u)$ in its neighborhood where $\PP'$-vertices connect to it, where~$n$ is a Poisson distributed random variable with parameter $\la' \abs{B(0, u)}$.
Finally, further $\PP$-vertices are generated by a Poisson process of intensity~$\la$ in the union of the neighborhoods $\bigcup_i B(y_i, v_i)$ of the $\PP'$-vertices.
Each of the neighborhoods $B(y_i, v_i)$ can be split to three parts, where a connecting $\PP$-vertex $(x, w)$ must satisfy the below conditions:
\[ \begin{array}{rll}
    \text{left tail}  & w \le (y_i - x) / \b v_i^{-\g'} & x < y_i - \b v_i^{-\g'}; \\
    \text{center}     & w \le 1                         & y_i - \b v_i^{-\g'} \le x \le y_i + \b v_i^{-\g'}; \\
    \text{right tail} & w \le (x - y_i) / \b v_i^{-\g'} & y_i + \b v_i^{-\g'} < x.
\end{array} \]
Next, we take the union of the center parts of the neighborhoods, and generate $\PP$-vertices in the union.
Finally, we determine the dominating tail parts at every position~$x$, and generate $\PP$-vertices in the corresponding tail parts.
Then, the simplices of the Palm distribution are those whose lowest mark $\PP$-vertex is $(0, u)$.

With the model parameters set, the simulation framework provides an opportunity to analyze the distribution of $\De_0$-degrees, $\De_1$-degrees, and $\De_0'$-degrees in the generated finite networks of different sizes by varying the intensities of the point processes.
Specifically, to examine the finite-size effect on the distributions of the network properties of interest, we simulated \num{100}~networks and kept the model parameters $\g = 0.7$, $\g' = 0.2$ constant, but altered the expected number~$\la$ and~$\la'$ of $\PP$- and $\PP'$-vertices, respectively.
To examine the Palm distribution, we simulated \num{100}~sets of \num{10000}~networks.
In this case, the values of $\la, \la'$ are indifferent, and we set them to $\la = \la' = 1$.
Lastly, the parameter~$\b$ was set so that the expected number of $\De_0$-degrees in the infinite network limit was~$3$.

As shown in \autoref{thm:deg}, all the examined distributions exhibit a power-law tail in the infinite network limit.
Thus, we fitted a power law distribution to the simulated values using the maximum likelihood method described in \citet{bauke2007parameter}.
Determining the minimum value $\xmin$ from which the power-law distribution should be fitted is not trivial.
If we use lower values of $\xmin$, more data points can be considered for the fitting, resulting in a more robust fit of the power-law distribution.
However, the estimated power-law exponent is more biased, as the power-law might not hold at lower values.
Conversely, for higher values of $\xmin$, the estimated power-law exponent has less bias, but the fit is less robust as fewer data points are considered.

\begin{figure} \centering
    \includegraphics{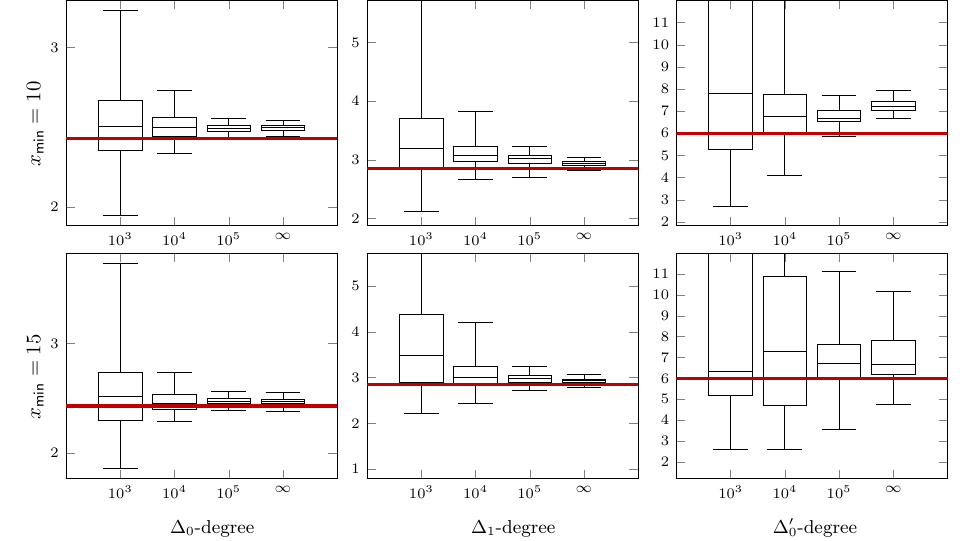}
    \vspace{-0.5cm}
    \caption{
        Fluctuation of the degree distribution exponents for different network sizes.
        The top row shows the power-law exponents when estimated from values larger than $\xmin = 10$.
        In contrast, the bottom row contains boxplots for the distribution of the exponents estimated from values larger than $\xmin = 15$.
    }\label{fig:degree_exponents}
\end{figure}

The boxplots of \autoref{fig:degree_exponents} show how the fitted power-law exponents of the above-mentioned degree distributions fluctuate relative to the theoretical value derived for infinite networks in \autoref{thm:deg}.
In the top row, we used $\xmin = 10$ for the fitting, while in the bottom row, we used $\xmin = 15$.
It can be observed that the exponents in the top row fluctuate less but exhibit a higher bias, whereas the exponents in the bottom row fluctuate more but display a lower bias.
This is because a set of $\PP$-vertices form a simplex if they share a common $\PP'$-vertex, and this more complex relationship between the $\PP$- and $\PP'$-vertices results in a higher variance of the estimated exponents.
We can also see that the degree exponents fluctuate much more for the $\De_0'$-degree distribution than for the $\De_0$-degree distribution.
This can be explained by the different decay rates of the neighborhoods of $\PP$- and $\PP'$-vertices: the power-law tails have exponents $1/\g = 1.43$ and $1/\g' = 0.5$, respectively, which results in fewer data points in the tail of the $\De_0'$-degree distribution compared to the $\De_0$-degree distribution.
As it will be shown in \autoref{sec:dat}, the $\De_0$- and $\De_0'$-degree distributions will be of particular interest when fitting the model parameters~$\g$ and~$\g'$ to the datasets.
Analyzing the boxplots of the $\De_0$- and $\De_0'$-degree distributions, we can see that the power-law exponents of the degree distributions are close to the theoretical values derived for infinite networks above a network size of \num{10000}.
Considering the optimal trade-off between bias and robustness, we chose to use $\xmin = 10$ for the fitting in the following sections.

\begin{figure} \centering
    \includegraphics{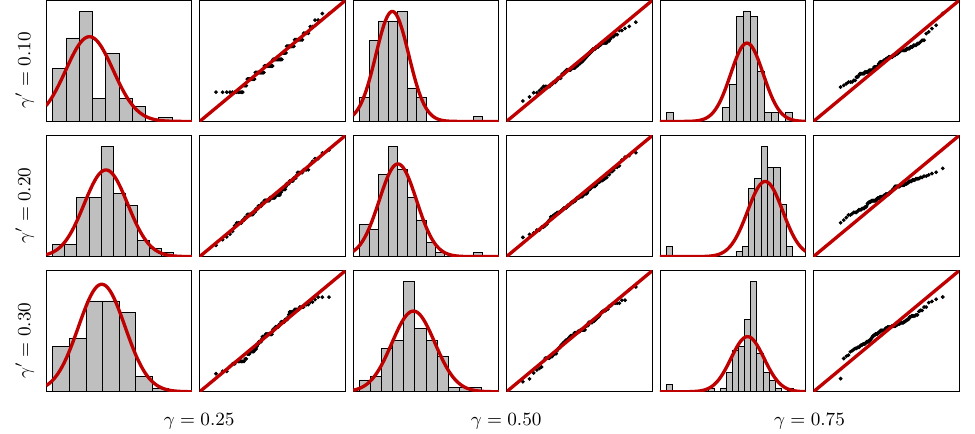}
    \caption{Distribution of the first Betti numbers with different~$\g$ and~$\g'$ parameters}\label{fig:betti_number_1_distributions}
\end{figure}

Next, we analyze the distributions of the first Betti number by simulating \num{100}~networks with fixed $\la = \la' = \num{100000}$ and varying~$\g$ and~$\g'$ parameters.
The histograms with fitted normal distributions and the corresponding Q-Q plots are shown in \autoref{fig:betti_number_1_distributions}.
In alignment with \autoref{thm:bet}, the Q-Q plots show that the Betti number is normally distributed for parameters $\g < 1/4$ and $\g' < 1/8$.
The figures for $1/4 \le \g \le 1/2$ and $\g' \ge 1/8$ suggest that the normal distribution of the first Betti number also holds for these parameter ranges.
The distribution of the Betti number for $\g = 0.75$, however, exhibits a fat left tail.
Note that the leftmost bins of the histograms for $\g = 0.75$ are at~$0$, which means that some simulated networks did not contain any loops.
This behavior can be explained by the presence of low-mark $\PP$- and/or $\PP'$-vertices that connect to many $\PP'$- and $\PP$-vertices when the edge count has an infinite variance, making the loops less likely to form.
This observation leads to the conjecture that the Betti numbers have $\a$-stable distributions if $\g > 0.5$.

\begin{figure} \centering
    \includegraphics{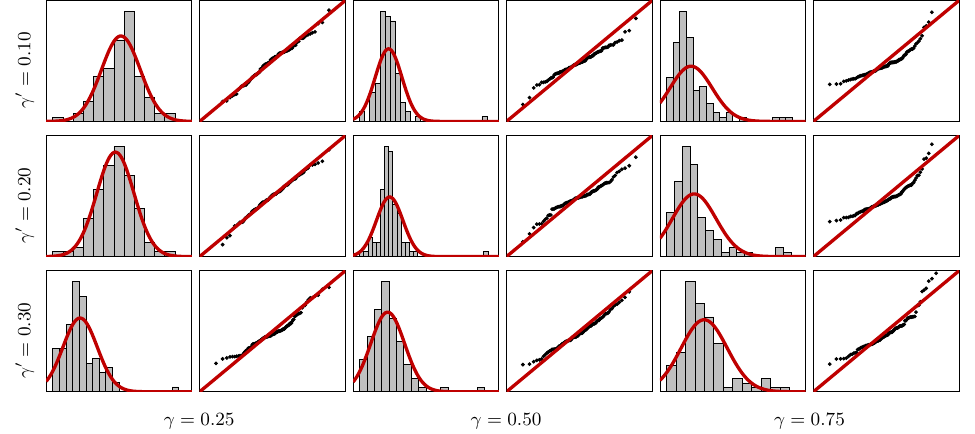}
    \caption{
        Distribution of edge counts with different~$\g$ and~$\g'$ parameters.
        For each combination of~$\g$ and~$\g'$, we show the distribution of edge counts with the fitted normal distribution and the corresponding Q-Q plot.
    }\label{fig:edge_count_distributions}
\end{figure}

\begin{figure} \centering
    \includegraphics{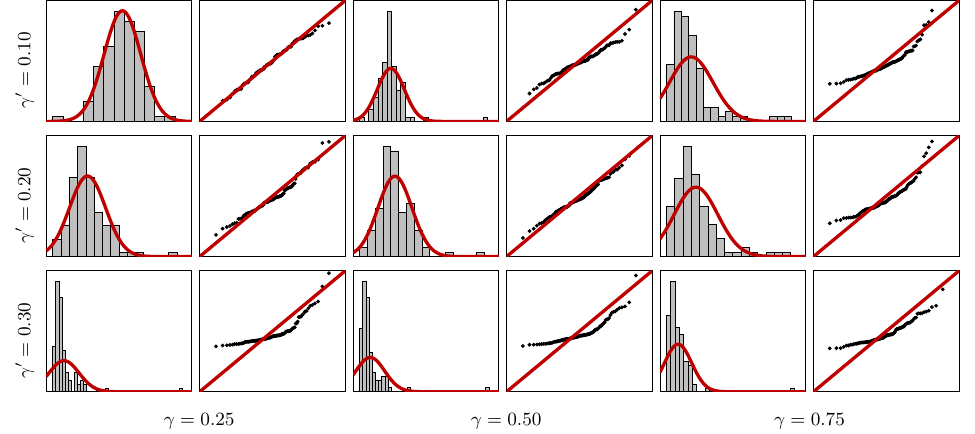}
    \caption{
        Distribution of triangle counts with various~$\g$ and~$\g'$ parameters.
        For each combination of~$\g$ and~$\g'$, we show the distribution of triangle counts with the fitted normal distribution and the corresponding Q-Q plot.}\label{fig:triangle_count_distributions}
\end{figure}

Finally, we analyze the distribution of the number of edges and triangles in the simulated networks.
As before, we simulated \num{100}~networks with fixed $\la = \la' = \num{100000}$ and varying~$\g$ and~$\g'$ parameters.
The distribution of the edge counts is presented in \autoref{fig:edge_count_distributions}, where we fitted a normal distribution to the simulated values and visualized the corresponding Q-Q plots.
As shown in \autoref{thm:stab} in the infinite network size limit, now the Q-Q plots show for finite networks that the number of edges are normally distributed for parameters $\g < 1/2$ and $\g' < 1/3$, and that the number of edges have an $\a$-stable distribution if $\g > 1/2$ and $\g' \ge 1/3$.
For $\g = 1/2$, it can be inferred from the Q-Q plots that the distribution of the number of edges is close to a normal distribution.
The distribution of the triangle counts together with the fitted normal distributions and the Q-Q plots are shown in \autoref{fig:triangle_count_distributions}.
Again, the simulated distributions of finite networks align with \autoref{thm:simp}: the Q-Q plots show that the number of triangles is normally distributed for parameters $\g < 1/2$ and $\g' < 1/5$.
It is also clear from the Q-Q plots of \autoref{fig:triangle_count_distributions} that if $\g > 1/2$ and $\g' < 1/5$, the number of triangles has an $\a$-stable distribution with $\a = 1/\g$, which was also shown for infinite networks in \autoref{thm:simp}.
We can furthermore conjecture that the number of triangles has an $\a$-stable distribution if $\g' > 1/5$ irrespective of the value of~$\g$.

After analyzing the degree distributions, Betti numbers, and the number of edges and triangles in the simulated networks, we can conclude that the theoretical results derived for infinite networks also hold for finite networks.
We finish the simulation study by analyzing further degree distributions, which, although not of particular interest in the following sections, may be examined in a further study.
\begin{itemize}
    \item{The $0$-coface degree distribution is the distribution of the number of edges incident to a $\PP$-vertex.}
    \item{The $1$-coface degree distribution is the distribution of the number of triangles incident to an edge.}
\end{itemize}
As we will see in \autoref{sec:dat}, these degree distributions exhibit a power-law tail in the datasets we analyze.

\begin{figure} \centering
    \includegraphics{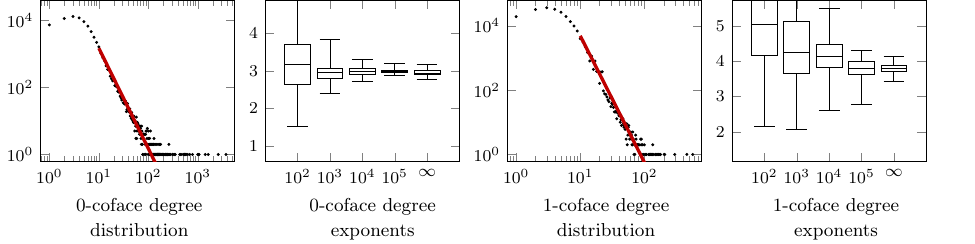}
    \vspace{-.5cm}
    \caption{
        Typical coface degree distributions and the fluctuation of the power-law exponents for different network sizes.
        The power-law exponents were estimated from values larger than $\xmin = 10$.
    }\label{fig:coface_degree_distributions}
\end{figure}

To examine the coface degree distribution of the networks generated by the model, we calculated the coface degree distributions for the sample network described at the beginning of \autoref{sec:sim}.
We also examined the fluctuations of the coface degree distribution exponents for different network sizes, and, as before, we simulated \num{100}~networks with varying $\la = \la'$ intensities and fixed $\g = 0.7$ and $\g' = 0.2$ parameters.
In this case again, we used $\xmin = 10$ for the fitting and set~$\b$ so that $\E{\De_0\ph{\ast}} = 3$.
The distributions of the $0$- and $1$-coface degrees and the fluctuations of the exponents are shown in \autoref{fig:coface_degree_distributions}.
We see that the distributions of the coface degrees are heavy-tailed, and the power-law exponents of the distributions converge above a network size of \num{10000}.

\section{Analysis of collaboration networks}\label{sec:dat}

After the analysis of the finite-size effects in \autoref{sec:sim}, we now compare our model with real-world datasets.
To do so, we chose to analyze coauthorship networks, where the authors are represented by $\PP$-vertices, and a $\PP'$-vertex represents each publication.

Using the Python wrapper for the arXiv API, we collected all documents from the arXiv preprint server that were uploaded up to 18 June 2024.
Based on the label of the documents referring to their primary categories, we created four datasets for the four largest scientific fields, namely, \texttt{computer science}, \texttt{electrical engineering and systems science}, \texttt{mathematics}, and \texttt{statistics}, that we denote by \texttt{cs}, \texttt{eess}, \texttt{math}, and \texttt{stat}, respectively.
In these datasets, we had access to the list of authors' names of the documents.
Identifying the authors by their names, we created
\begin{itemize}
    \item{a bipartite network where the $\PP$-vertices are the authors, and the $\PP'$-vertices are the documents;}
    \item{a simplicial complex representing the higher-order network of the authors.}
\end{itemize}
As simplicial complexes are closed under taking subcomplexes, each document with~$n$ authors is represented by $2^n - 1$-simplices in the simplicial complex.
This exponential dependence of the number of simplices on the number of authors led to a disproportionate influence of documents with a high number of authors on our results.
To mitigate this effect, we included documents with at most~$20$ authors in our analysis, thereby neglecting a proportion of $0.13\%$ of the documents.

The main properties of the datasets are summarized in \autoref{tab:dataset_properties}.
\begin{table} \centering
    \sisetup{group-separator = {}, group-minimum-digits = 4}
    \caption{Main properties of the datasets} \label{tab:dataset_properties}
    \begin{tabular}{l S[table-format=6.0] S[table-format=6.0] S[table-format=5.0] S[table-format=6.0]} \toprule
        dataset & {authors} & {documents} & {components} & {size of largest component} \\ \midrule
        \texttt{cs}      &    504959 &      540573 &        23909 &                      439375 \\
        \texttt{eess}    &     93005 &       85445 &         6061 &                       67700 \\
        \texttt{math}    &    211757 &      500403 &        27157 &                      163738 \\
        \texttt{stat}    &     48700 &       40869 &         4236 &                       36094 \\
        \bottomrule
    \end{tabular}
\end{table}

\begin{figure} \centering
    \includegraphics{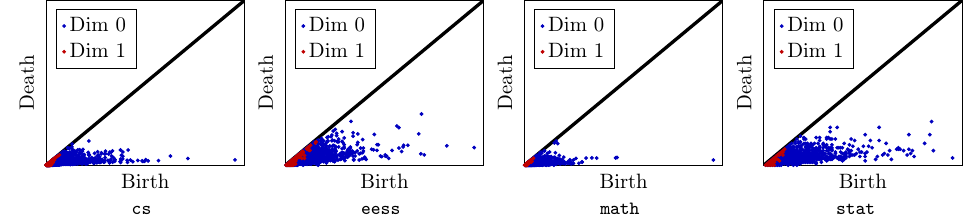}
    \caption{Persistence diagrams of the datasets.}\label{fig:persistence_diagrams}
\end{figure}

Using the documents as $\PP'$-vertices, we created a filtered simplicial complex for each dataset.
We set the filtration value of a set of authors to be the number of papers they are all authors of, and used a filtration value decreasing from~$\ff$ to~$0$ to create the filtered Dowker complex.
We present the persistence diagrams computed from the filtered Dowker complex in \autoref{fig:persistence_diagrams}.
Although the datasets contain many components and loops, some of them share the same birth and death filtration values, resulting in fewer points in the persistence diagrams.

\begin{figure}[p] \centering
    \includegraphics{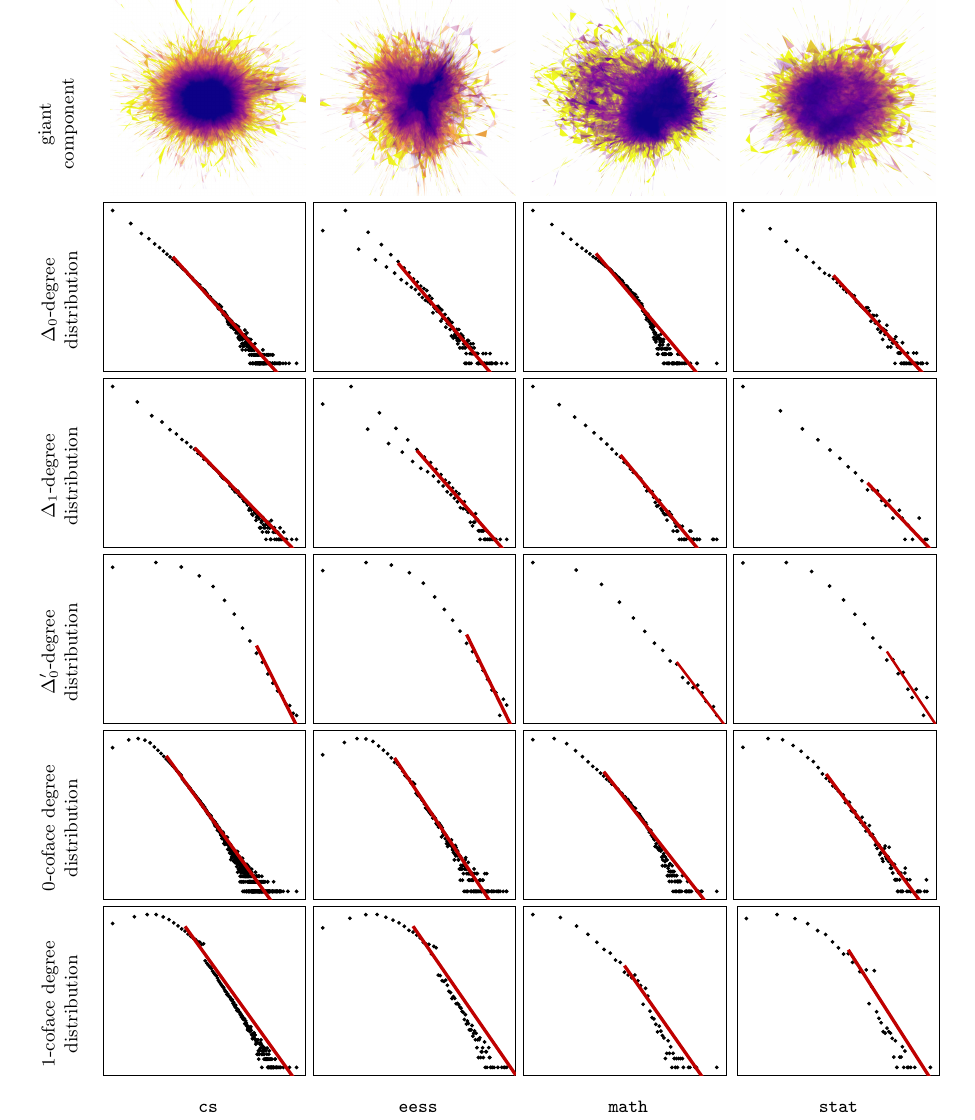}
    \caption{Main properties of the datasets}\label{fig:dataset_main_properties}
\end{figure}

Further properties of the datasets are illustrated in \autoref{fig:dataset_main_properties}.
In the top row, we visualized the largest components of the datasets, where each document is represented by a polygon around the authors, with darker colors indicating documents with more authors.
Further rows display the degree distributions of the datasets, along with the fitted power-law distributions.
The degree distributions that we examine are similar to those before in \autoref{sec:sim}:
\begin{itemize}
    \item{The $\De_0$-degree distribution is the distribution of the number of documents an author published.}
    \item{The $\De_1$-degree distribution describes the distribution of the number of documents that pairs of authors wrote together.}
    \item{The distribution of the number of authors of a document is the $\De_0'$-degree distribution.}
\end{itemize}
We observe that all the visualized degree distributions exhibit heavy tails, and the fitted power-law distributions accurately describe the data.

To compare the datasets with our model, we first need to estimate the model parameters.
The~$\g$ and~$\g'$ parameters are estimated from the exponents of the fitted power-law distributions to the $\De_0$- and $\De_0'$-degree distributions, respectively.

Note that the datasets do not contain isolated $\PP$- and $\PP'$-vertices, i.e., authors without any documents and documents without any authors.
Thus, we compensate for these in our model so that the number of authors and documents in the datasets match the number of non-isolated $\PP$- and $\PP'$-vertices in the simulations, respectively.

Let $\PP_0$ and $\PP'_0$ be the set of isolated $\PP$- and $\PP'$-vertices in the simulated simplicial complex, respectively, and let $p \coeq (x, u) \in \PP$ and $p' \coeq (z, w) \in \PP'$.
Then, the expected number of non-isolated $\PP$-vertices in the simulation is given by
\[ \begin{aligned}
    \E[\big]{\#(\PP \sm \PP_0)} &= \E[\Big]{\sum_{p \in \PP} \1[\big]{\PP'(B(p)) \ne 0}} = \la \int_0^1 \E[\big]{\1[\big]{\PP'(B(p)) \ne 0}} \d u \\
    &= \la \int_0^1 1 - \P{\PP'(B(0, u)) = 0} \d u = \la \int_0^1 1 - \exp[\big]{- \la' \abs{B(0, u)}} \d u \\
    &= \la \int_0^1 1 - \exp[\bigg]{- \la' \int_0^1 \int_{- \b u^{-\g} w^{-\g'}}^{\b u^{-\g} w^{-\g'}} \d z \d w} \d u = \la \biggl( 1 - \int_0^1 \exp[\Big]{- \f{2 \b \la'}{1 - \g'} u^{-\g}} \d u \biggr).
\end{aligned} \]
After substituting $\mu \coeq 2 \b \la' / (1 - \g') u^{-\g}$, we obtain
\[ \begin{aligned}
    \E[\big]{\#(\PP \sm \PP_0)} &= \la \Bigl( 1 - \f{1}{\g} \Bigl( \f{1 - \g'}{2 \b \la'} \Bigr)^{- \f{1}{\g}} \int_{\f{2 \b \la'}{1 - \g'}}^\ff \mu^{- \f{1}{\g} - 1} \e^{- \mu} \d \mu \Bigr) \\
    &= \la \Bigl( 1 - \f{1}{\g} \Bigl( \f{1 - \g'}{2 \b \la'} \Bigr)^{- \f{1}{\g}} \Ga \Bigl( - \f{1}{\g}, \f{2 \b \la'}{1 - \g'} \Bigr) \Bigr),
\end{aligned} \]
where~$\Ga$ denotes the upper incomplete gamma function.
Following the same steps, we find that the expected number of non-isolated $\PP'$-vertices in the simulation is given by
\[ \E[\big]{\#(\PP' \sm \PP'_0)} = \la' \Bigl( 1 - \f{1}{\g'} \Bigl( \f{1 - \g}{2 \b \la} \Bigr)^{- \f{1}{\g'}} \Ga \Bigl( - \f{1}{\g'}, \f{2 \b \la}{1 - \g} \Bigr) \Bigr). \]
Furthermore, we set the parameter~$\b$ so that the expected number of edges in~$\Gb$ in the simulation matches the number of author--document edges in the datasets.
Additionally, by the Mecke formula, we have
\[ \E[\Big]{\sum_{\substack{p \in \PP \\ p' \in \PP'}} \1[\big]{p \in B(p')}} = \int_{\S^2} \P{p \in B(p')} \d p \d p' = \f{2 \b \la \la'}{(1 - \g)(1 - \g')}. \]
By numerically solving the above system of equations given by the expected number of non-isolated $\PP$- and $\PP'$-vertices, we obtain the parameter estimates shown in \autoref{tab:dataset_parameter_estimates}.
\begin{table} \centering
    \caption{Fitted power-law exponents and the inferred model parameters} \label{tab:dataset_parameter_estimates}
    \sisetup{group-separator = {}, group-minimum-digits = 4}
    \begin{tabular}{l S[table-format=1.2] S[table-format=1.2] S[table-format=1.2] S[table-format=1.2] S[table-format=1.2e1] S[table-format=6.0] S[table-format=6.0]} \toprule
        \multirow{2}{*}{dataset \hspace{0.3cm}} & \multicolumn{2}{c}{{$\De_0$-degree}} & \multicolumn{2}{c}{{$\PP'$-vertex degree}} & \multicolumn{1}{c}{\multirow{2}{*}{{$\b$}}} & \multicolumn{1}{c}{\multirow{2}{*}{{$\la$}}} & \multicolumn{1}{c}{\multirow{2}{*}{{$\la'$}}} \\ \cmidrule(lr){2-3} \cmidrule(lr){4-5}
                & {exponent} & {$\g$} & {exponent} & {$\g'$} &         &        &        \\ \midrule
        \texttt{cs}      &       2.37 &   0.73 &       5.46 &    0.22 & 8.19e-7 & 579719 & 532491 \\
        \texttt{eess}    &       2.75 &   0.57 &       5.50 &    0.22 & 7.89e-6 &  98528 &  83985 \\
        \texttt{math}    &       2.44 &   0.69 &       6.51 &    0.18 & 1.03e-6 & 231606 & 588628 \\
        \texttt{stat}    &       2.89 &   0.53 &       5.74 &    0.21 & 1.13e-5 &  57488 &  41655 \\
        \bottomrule
    \end{tabular}
\end{table}

\begin{figure} \centering
    \includegraphics{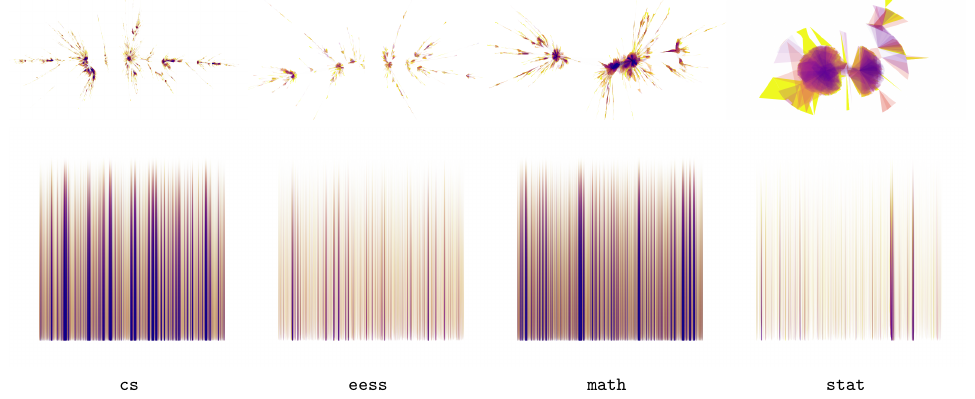}
    \caption[Simplicial complexes built from the datasets]{The largest components of the simplicial complexes built from the datasets}\label{fig:model_sample_plots}
\end{figure}

Using these parameters, we simulated a sample network for each dataset, which is visualized in \autoref{fig:model_sample_plots}.
While the top row contains the largest components of the datasets, the bottom row displays the $\PP$-vertices in the space~$\S$.
In both rows, as above, each $\PP'$-vertex is represented by a polygon around $\PP$-vertices connecting to them with a darker color if the $\PP'$-vertex is connected to more $\PP$-vertices.
We can see that the model networks possess a tree-like structure, dominated by a few high-degree $\PP$-vertices.

To determine if our model can capture the higher-order structures of the datasets, we need to perform hypothesis tests on the simplex counts and the first Betti number.
Although \autoref{thm:stab} and our conjecture in \autoref{sec:sim} show that $\a$-stable distributions describe the simplex counts and the first Betti number, it is not clear how to fully specify all parameters of these distributions.
Therefore, we estimated the parameters of the $\a$-stable distributions from the simulated networks to statistically compare the datasets with those generated by the model.
For each dataset, we simulated~$100$ networks, and fitted an $\a$-stable distribution to the values with $\a = 1 / \g$, which described the data well in all the cases.

\begin{figure} \centering
    \includegraphics{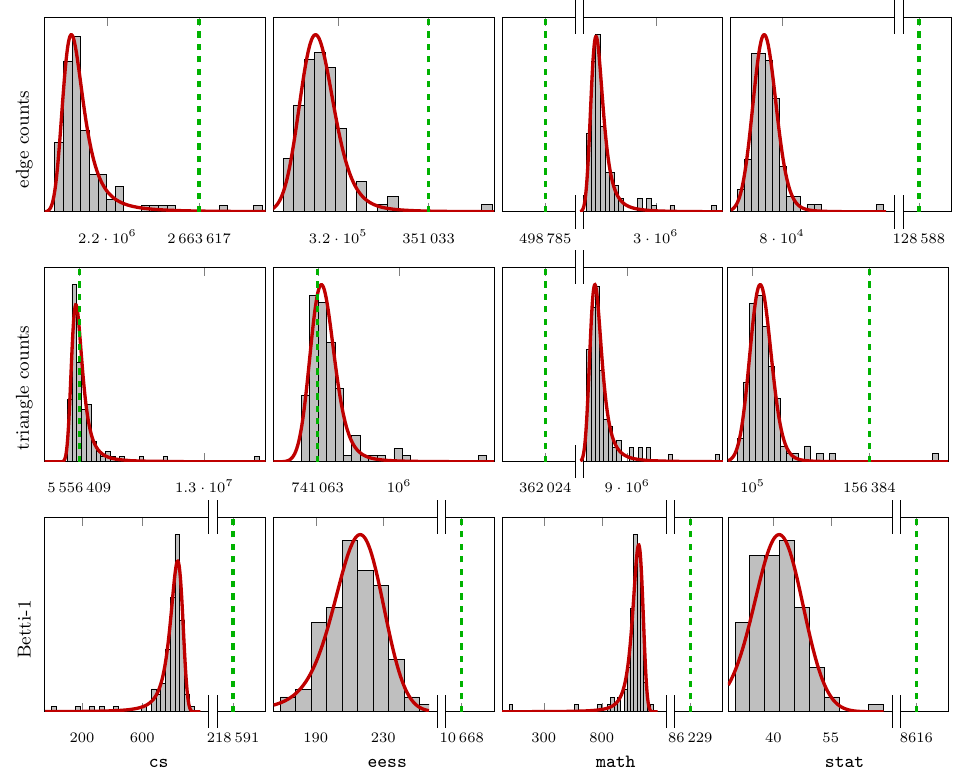}
    \vspace{-.5cm}
    \caption[Hypothesis testing of the edge and triangle counts and Betti-1]{
        Hypothesis testing of the edge counts (top), triangle counts (center), and Betti-1 (bottom) for the datasets.
        The model parameters were determined based on the parameters of the datasets described in \autoref{tab:dataset_parameter_estimates}.
    }\label{fig:simplex_count_hypothesis_test}
    \vspace{-.2cm}
\end{figure}

The results of the hypothesis tests for the simplex counts and the first Betti numbers are summarized in \autoref{fig:simplex_count_hypothesis_test}, and the parameters of the fitted distributions are shown in Tables~\ref{tab:num_of_edges_hypothesis_tests}, \ref{tab:num_of_triangles_hypothesis_tests}, and~\ref{tab:betti_number_1_hypothesis_tests}.
Although the model parameters were set to ensure that the expected number of $\De_0$-degrees in the simulation matches the number of author-document edges in the datasets, the hypothesis tests indicate that the higher-order structures of the datasets differ from those of the generated networks.
The case of the \texttt{eess} and \texttt{math} datasets is particularly interesting, as in the case of the former, the simulated networks contain more triangles, but fewer edges, compared to the datasets.
In contrast, in the case of the latter, the opposite is true.
In the case of the first Betti numbers, the values of the datasets are substantially larger than for the generated networks, and all the hypothesis tests are rejected.
The reason for this deviation is that the simulated networks are dominated by a few high-degree $\PP$-vertices, resulting in a tree-like structure with a relatively small number of loops.
At the same time, the datasets contain more complex structures.
\begin{table} \centering
    \caption{Results of the hypothesis tests for the edge counts}
    \label{tab:num_of_edges_hypothesis_tests}
    \sisetup{group-separator = {}, group-minimum-digits = 4}
    \begin{tabular}{l S[table-format=7.0] S[table-format=1.4] S[table-format=1.1] S[table-format=7.0] S[table-format=5.0] S[table-format=1.2e1]} \toprule
        dataset & {dataset value} & {$\hat{\a}$} & {$\hat{\b}$} & {location} & {scale} & {$p$-value} \\ \midrule
        \texttt{cs}      &         2663617 &       1.3658 &          1.0 &    2093994 &   43104 &    2.69e-02 \\
        \texttt{eess}    &          351033 &       1.7506 &          1.0 &     314266 &    4331 &    1.13e-02 \\
        \texttt{math}    &          498785 &       1.4415 &          1.0 &    2601637 &   36355 &    0.00e+00 \\
        \texttt{stat}    &          128588 &       1.8918 &          1.0 &      77681 &    1201 &    1.73e-04 \\
        \bottomrule
    \end{tabular}
\end{table}

\begin{table} \centering
    \caption{Results of the hypothesis tests for the number of triangles}
    \label{tab:num_of_triangles_hypothesis_tests}
    \sisetup{group-separator = {}, group-minimum-digits = 4}
    \begin{tabular}{l S[table-format=7.0] S[table-format=1.4] S[table-format=1.1] S[table-format=7.0] S[table-format=6.0] S[table-format=1.2e1]} \toprule
        dataset & {dataset value} & {$\hat{\a}$} & {$\hat{\b}$} & {location} & {scale} & {$p$-value} \\ \midrule
        \texttt{cs}      &         5556409 &       1.3658 &          1.0 &    5823007 &  263597 &     9.13e-1 \\
        \texttt{eess}    &          741063 &       1.7506 &          1.0 &     766277 &   28281 &     6.50e-1 \\
        \texttt{math}    &          362024 &       1.4415 &          1.0 &    6973583 &  187554 &     0.00e+0 \\
        \texttt{stat}    &          156384 &       1.8918 &          1.0 &     104413 &    3720 &     1.45e-3 \\
        \bottomrule
    \end{tabular}
\end{table}

\begin{table} \centering
    \sisetup{group-separator = {}, group-minimum-digits = 4}
    \caption{Results of the hypothesis tests for the first Betti numbers}
    \label{tab:betti_number_1_hypothesis_tests}
    \begin{tabular}{l S[table-format=6.0] S[table-format=1.4] S[table-format=2.1] S[table-format=4.0] S[table-format=2.0] S[table-format=1.2]} \toprule
        dataset & {dataset value} & {$\hat{\a}$} & {$\hat{\b}$} & {location} & {scale} & {$p$-value} \\ \midrule
        \texttt{cs}      &          218591 &       1.3658 &         -1.0 &        777 &      33 &        0.00 \\
        \texttt{eess}    &           10668 &       1.7506 &         -1.0 &        211 &      11 &        0.00 \\
        \texttt{math}    &           86229 &       1.4415 &         -1.0 &       1069 &      35 &        0.00 \\
        \texttt{stat}    &            8616 &       1.8918 &         -1.0 &         41 &       5 &        0.00 \\ \bottomrule
    \end{tabular}
\end{table}

The results of the hypothesis tests indicate that, although our model captures many properties of the datasets, it is too simplistic to accurately describe the higher-order structures of the datasets.

\FloatBarrier

\subsection*{Disclosure statement}

No potential conflict of interest was reported by the authors.

\subsection*{Funding}

This work was supported by the Danish Data Science Academy, which is funded by the Novo Nordisk Foundation (NNF21SA0069429) and Villum Fonden (40516).

\bibliographystyle{abbrv}

\end{document}